\documentclass[journal]{IEEEtran}

\ifCLASSINFOpdf
\else
\fi
\usepackage{amsmath,amssymb}
\usepackage{xcolor}
\usepackage{graphicx}
\usepackage{epstopdf}

\usepackage{etoolbox}
\makeatletter
\patchcmd{\@makecaption}
  {\scshape}
  {}
  {}
  {}

\newcommand{\thickhline}{%
    \noalign {\ifnum 0=`}\fi \hrule height 1pt
    \futurelet \reserved@a \@xhline
}
\makeatother

\def\red#1{\textcolor{black}{#1}}
\def\sa#1{\red{#1}}

\def\magenta#1{\textcolor{black}{#1}}
\def\grad{\nabla}

\def\cA{\mathcal{A}}

\def\cG{\mathcal{G}}

\def\cL{\mathcal{L}}

\def\cN{\mathcal{N}}
\def\cO{\mathcal{O}}

\def\cU{\mathcal{U}}

\def\cX{\mathcal{X}}

\def\mE{\mathbb{E}}

\def\mR{\mathbb{R}}

\def\smskip{\smallskip}

\def\texitem#1{\par\smskip\noindent\hangindent 25pt
               \hbox to 25pt {\hss #1 ~}\ignorespaces}

\def\norm#1{\left\|#1\right\|}

\newcommand{\BEAS}{\begin{eqnarray*}}
\newcommand{\EEAS}{\end{eqnarray*}}
\newcommand{\BEA}{\begin{eqnarray}}
\newcommand{\EEA}{\end{eqnarray}}
\newcommand{\BEQ}{\begin{eqnarray}}
\newcommand{\EEQ}{\end{eqnarray}}
\newcommand{\BIT}{\begin{itemize}}
\newcommand{\EIT}{\end{itemize}}
\newcommand{\BNUM}{\begin{enumerate}}
\newcommand{\ENUM}{\end{enumerate}}

\newcommand{\BA}{\begin{array}}
\newcommand{\EA}{\end{array}}
\newcommand{\BE}{\begin{equation}}
\newcommand{\EE}{\end{equation}}

\newcommand{\reals}{\mathbb{R}}
\newcommand{\integers}{\mathbb{Z}}

\newcommand{\Rank}{\mathop{\bf rank}}

\newcommand{\diag}{\mathop{\bf diag}}

\newif\ifpagenumbering
\pagenumberingtrue

\pagenumberingfalse

\newsavebox{\theorembox}
\newsavebox{\lemmabox}
\newsavebox{\defnbox}
\newsavebox{\corollarybox}
\newsavebox{\propositionbox}
\newsavebox{\remarkbox}
\newsavebox{\assbox}
\savebox{\theorembox}{\noindent\bf Theorem}
\savebox{\lemmabox}{\noindent\bf Lemma}
\savebox{\defnbox}{\noindent\bf Definition}
\savebox{\corollarybox}{\noindent\bf Corollary}
\savebox{\propositionbox}{\noindent\bf Proposition}
\savebox{\remarkbox}{\noindent\bf Remark}
\savebox{\assbox}{\noindent\bf Assumption}

\newtheorem{remark}{\usebox{\remarkbox}}[section]

\newtheorem{defn}{\usebox{\defnbox}}

\def\fprod#1{\langle#1 \rangle}
\DeclareMathOperator*{\argmin}{\arg\!\min}

\newcommand{\rank}{\mathbf{rank}}
\newcommand{\st}{\textrm{s.t.}}
\newcommand{\etal}{et al. }
\newcommand{\Ind}{\mathbf{1}}
\newcommand{\half}{\frac{1}{2}}
\newcommand{\sign}{\textrm{sign}}
\newcommand{\St}{{\mathcal{S}_t}}

\graphicspath{{./figures/}}

\begin{document}
\title{Efficient Optimization Algorithms for Robust Principal Component Analysis and Its Variants}

\author{Shiqian~Ma,
        and~Necdet~Serhat~Aybat
\thanks{S. Ma is with the Department of Mathematics at University of California, Davis, CA, 95616 USA e-mail: (sqma@math.ucdavis.edu).}
\thanks{N. S. Aybat is with the Department of Industrial and Manufacturing Engineering at Penn State, University Park, PA, 16802 USA e-mail: (nsa10@psu.edu).}
\thanks{Manuscript received January xx, 2018; revised xxx xx, 2018.}}

\markboth{Manuscript,~Vol.~xx, No.~xx, xxx 2018}
{Shell \MakeLowercase{\textit{et al.}}: Efficient Optimization Algorithms for Robust PCA}

\maketitle

\begin{abstract}
Robust PCA has drawn significant attention in the last decade due to its success in numerous application
domains, ranging from bio-informatics, statistics, and machine learning to image and video processing in computer
vision. Robust PCA and its variants such as sparse PCA and stable PCA can be formulated as optimization problems
with {exploitable} special structures. Many specialized efficient optimization methods have been proposed to solve robust PCA and
related problems. In this {paper} we review existing optimization methods for solving convex and nonconvex relaxations/variants of robust PCA, discuss their advantages and disadvantages, and elaborate on their convergence behaviors. We also provide some insights for possible future research directions including new algorithmic frameworks that might be suitable for implementing on multi-processor setting to handle large-scale problems.
\end{abstract}


\IEEEpeerreviewmaketitle

\section{Introduction}

\IEEEPARstart{P}{rincipal} component analysis (PCA) is a fundamental tool in statistics and data science. It obtains a low-dimensional expression for high-dimensional data in an $\ell_2$ sense. However, it is known that the classical PCA is sensitive to gross errors. Robust PCA (RPCA) has been proposed to remove the effect of sparse gross errors. For a given data matrix $M\in\mR^{m\times n}$, RPCA seeks to decompose it \sa{into} two parts $M:=\sa{L^o+S^o}$ \sa{where} $\sa{L^o}$ \sa{is} a low-rank matrix and $\sa{S^o}$ \sa{is} a sparse matrix. That is, RPCA assumes that $M$ is a superposition of $\sa{L^o}$ and $\sa{S^o}$. As a result, the gross errors will be captured by the sparse matrix $S^o$ so that the low-rank matrix $L^o$ can still approximate $M$ well. RPCA does not only provide a low-dimensional approximation which is robust to outliers, it also finds vast applications in a variety of real applications such as computer vision \cite{Candes-Li-Ma-Wright-RPCA}, \sa{image alignment} \cite{Peng-rasl-2012}, subspace recovery \cite{Liu-robust-subspace-recovery-2013}, clustering \cite{Shahid-iccv} and so on.

\sa{Mathematically, \cite{Candes-Li-Ma-Wright-RPCA,Chandrasekaran-RPCA-2011,NIPS2009_3704,NIPS2010_4005} investigated the conditions on the low-rank and sparse components $({L^o},{S^o})$ so that the inverse problem of recovering unknown $({L^o},{S^o})$ given $M$ is well defined. One particular formulation of RPCA can be stated as follows:}
\BE\label{rpca-rank-L0}
\min_{L,S\in\mR^{m\times n}} \rank(L) + \rho\|S\|_0 \ \st\ L+S = M,
\EE
where $\|S\|_0$ is called the $\ell_0$-norm\footnote{Technically, it is not a norm \sa{because it is not homogeneous;} but, we still call it a norm following the convention.} of $S$ and counts the number of nonzero entries of $S$, and $\rho>0$ is a 
\sa{tradeoff} parameter. It is known that \eqref{rpca-rank-L0} is NP-hard and thus numerically intractable. 
{Later,} it was shown in \cite{Candes-Li-Ma-Wright-RPCA,Chandrasekaran-RPCA-2011,NIPS2009_3704} that under certain conditions, \eqref{rpca-rank-L0} is equivalent to the following convex program with high probability:
\BE\label{rpca-nuclear-L1}
\min_{L,S\in\mR^{m\times n}} \|L\|_* + \rho\|S\|_1 \ \st\ L+S = M,
\EE
where $\|L\|_*$ is called the nuclear norm of $L$ and equal to the sum of the singular values of $L$, and $\|S\|_1:=\sum_{ij}|S_{ij}|$ is called the $\ell_1$ norm of $S$. \sa{The optimization problem in}~\eqref{rpca-nuclear-L1} is called robust principal component pursuit (RPCP), and it can be reformulated as a semidefinite program (SDP)~\sa{\cite{recht2010guaranteed}} and solved by \sa{an} interior point method \sa{for SDPs}. However, RPCA problems arising 
{in} practice are usually of very large scale, and interior point methods do not scale well for these problems. More efficient algorithms that solve \eqref{rpca-nuclear-L1} and its variants {by exploiting the structure} in these problems were studied extensively in the literature. One variant of \eqref{rpca-nuclear-L1} 
\sa{deals with an additional dense noise component. In particular, when} 
\sa{$M$ contains also a} dense noise component {$N^o$ such that $\norm{N^o}_F\leq \sigma$ for some noise level $\sigma>0$}, i.e., $M=L^o+S^o+N^o$, \sa{instead of RPCP formulation} in \eqref{rpca-nuclear-L1}, 
the following so-called stable PCP (SPCP) problem \sa{is solved:}
\BE\label{rpca-nuclear-L1-noise}
\min_{L,S\in\mR^{m\times n}} \|L\|_* + \rho\|S\|_1 \ \st\ \norm{L+S-M}_F\leq\sigma.
\EE
{It is proved in \cite{Zhou-stable-pca-2010} that, {under certain conditions on $M$,} solving \eqref{rpca-nuclear-L1-noise} gives a stable estimate of $L^o$ and $S^o$ with high probability in the sense that $\|\hat{L}-L^o\|_F^2 + \|\hat{S}-S^o\|_F^2 \leq {\cO(mn\sigma^2)}$
where $(\hat{L},\hat{S})$ denotes the optimal solution to \eqref{rpca-nuclear-L1-noise}.}
\sa{Since \eqref{rpca-nuclear-L1-noise} satisfies the Slater's condition, it is equivalent to the following unconstrained problem for an appropriately chosen penalty parameter {$\mu>0$} depending on $\sigma$: 
\BE\label{rpca-nuclear-L1-noise-reform-unconstrained}
\min_{L,S\in\mR^{m\times n}} \|L\|_* + \rho\|S\|_1 + \frac{{\mu}}{2} \|L+S-M\|_F^2.
\EE}

Note that if $M$ is only partially observed, that is, if we only have observations on $M_{ij}$ for some indices $(i,j)$ \sa{from} a subset $\Omega$, then \eqref{rpca-nuclear-L1}, \eqref{rpca-nuclear-L1-noise} and \eqref{rpca-nuclear-L1-noise-reform-unconstrained} can be respectively \sa{reformulated as} 
\BE\label{rpca-nuclear-L1-partial}
\min_{L,S\in\mR^{m\times n}} \|L\|_* + \rho\|S\|_1 \ \st\ \sa{P_\Omega(L+S-M) = \mathbf{0},}
\EE
\BE\label{rpca-nuclear-L1-noise-partial}
\min_{L,S\in\mR^{m\times n}} \|L\|_* + \rho\|S\|_1 \ \st\ \sa{\norm{P_\Omega(L+S-M)}_F \leq \sigma,}
\EE
\BE\label{rpca-nuclear-L1-noise-reform-unconstrained-partial}
\min_{L,S\in\mR^{m\times n}} \|L\|_* + \rho\|S\|_1 + \frac{{\mu}}{2} \|P_\Omega(L+S-M)\|_F^2,
\EE
where the operator $P_\Omega: \mR^{m\times n}\rightarrow \mR^{m\times n}$ is defined as $[P_\Omega(M)]_{ij} = M_{ij}$, if $(i,j)\in\Omega$, and $[P_\Omega(M)]_{ij}=0$ otherwise. Most algorithms we discuss in this {paper} for solving \eqref{rpca-nuclear-L1}, \eqref{rpca-nuclear-L1-noise} and \eqref{rpca-nuclear-L1-noise-reform-unconstrained} can be used to solve \eqref{rpca-nuclear-L1-partial}, \eqref{rpca-nuclear-L1-noise-partial} and \eqref{rpca-nuclear-L1-noise-reform-unconstrained-partial} directly or with very little modification. For brevity, we will only describe algorithms for solving \eqref{rpca-nuclear-L1}, \eqref{rpca-nuclear-L1-noise} and \eqref{rpca-nuclear-L1-noise-reform-unconstrained} in this {paper}.

{For the sake of completeness,} \sa{we here briefly describe the results in \cite{Candes-Li-Ma-Wright-RPCA,Zhou-stable-pca-2010}. Assume that given data matrix $M\in\reals^{m\times n}$ is a superposition of the \emph{unknown} component matrices $L^o$, $S^o$, and $N^o$, i.e., $M=L^o+S^o+N^o$, such that $L^o$ is low-rank ($r:=\rank(L^0)\ll\min\{m,n\}$), $S^o$ is sparse ($s:=\norm{S^o}_0\ll mn$), and $\norm{N^o}\leq\sigma$ for some $\sigma>0$. Robust/stable PCA is an inverse problem with the objective of recovering $L^o$ and $S^o$ from the data matrix $M$. Clearly, if $S^o$ is low-rank and/or $L^o$ is sparse, the recovery is hopeless. To avoid these pathological instances, \cite{Candes-Li-Ma-Wright-RPCA,Zhou-stable-pca-2010} consider $(L^o,S^o)$ pairs coming from a particular class satisfying some incoherence and randomness conditions. Suppose the singular value decomposition (SVD) of $L^o$ is given by $L^o = \sum_{i=1}^r \sigma_iu_iv_i^\top$, where
$U=[u_1,\ldots,u_r]$ and $V = [v_1,\ldots,v_r]$ are formed by its left- and right-singular vectors. The incoherence conditions assume that there exists a parameter {$\delta>0$} such that
\begin{subequations}\label{incoherence-conditions}
\begin{align}
\max_i\|U^\top e_i\|^2\leq {\delta} r/m,\ & \max_i\|V^\top e_i\|^2 \leq {\delta} r/n, \label{incoherence-conditions-1} \\
\|UV^\top\|_\infty\leq\sqrt{\frac{{\delta} r}{mn}},\  & \label{incoherence-conditions-2}
\end{align}
\end{subequations}
where $\|Z\|_\infty := \max_{ij}|Z_{ij}|$ and $e_i$ denotes the $i$-th unit vector. For $S^o$, it is assumed that the set of indices for the nonzero entries is random and follows a uniform distribution among all the subsets of cardinality $s$. Roughly speaking, these conditions assure that the low-rank matrix $L^o$ is not sparse, and the sparse matrix $S^o$ is not low-rank.
Under these assumptions
it is shown in~\cite{Candes-Li-Ma-Wright-RPCA} that when $N^o=\mathbf{0}$, i.e., $\sigma=0$, solving the convex program \eqref{rpca-nuclear-L1} with $\rho=\frac{1}{\sqrt{\max\{m,n\}}}$ recovers the optimal solution of \eqref{rpca-rank-L0}, $(L^*,S^*)$, with high probability and $(L^*,S^*)=(L^o,S^o)$ provided that $L^o$ is sufficiently low-rank and $S^o$ is sufficiently sparse comparing to the matrix size -- see also~\cite{NIPS2009_3704}.} \sa{These results extend to the case where $M$ is partially observed; indeed, Cand\`{e}s et al.~\cite{Candes-Li-Ma-Wright-RPCA} show that solving~\eqref{rpca-nuclear-L1-partial} recovers $L^o$ under similar conditions.} Moreover, in \cite{Zhou-stable-pca-2010}, the authors showed that under certain conditions, solving the convex problem \eqref{rpca-nuclear-L1-noise} again with $\rho=\frac{1}{\sqrt{\max\{m,n\}}}$ \sa{generates a low-rank and sparse decomposition $(L^*,S^*)$ such that $\norm{L^*-L^o}_F^2+\norm{S^*-S^o}_F^2\leq Cmn{\sigma}^2$ for some {constant $C>0$ (independent of $m$, $n$ and $\sigma$)} with high probability -- note when $N=\mathbf{0}$, the recovery is exact with high probability.} 

\red{Recently, there are works that further study statistical guarantees of different RPCA models. Zhang, Zhou and Liang \cite{Zhang-zhou-Liang} provide a refined analysis of RPCA which allows the support of the \sa{sparse} error matrix to be generated with non-uniform sampling, i.e., entries of the low-rank matrix are corrupted 
with 
different probabilities -- \sa{hence, one can model the scenario where some entries are more prone to corruption than the others}. A nonconvex model of RPCA is studied in \cite{Yi-Park-Chen-Caramanis-gradient} and a gradient descent method with proper initialization is shown to be able to reduce the computational complexity comparing with existing methods. Zhang and Yang \cite{Zhang-Yang-manifold-2017} consider a nonconvex optimization formulation with manifold constraint for RPCA. Two algorithms for manifold optimization are proposed in \cite{Zhang-Yang-manifold-2017} and it is shown that they can reduce the dependence on the \sa{condition} number of the underlying low-rank matrix theoretically. Netrapalli \etal \cite{Netrapalli-nonconvex-2014} consider another nonconvex formulation of RPCA and analyze the iteration complexity of the proposed alternating projection method.}

{There are also recent survey papers \cite{Bouwmans-survey-cviu-2014,Aybat-book-2016} that discuss algorithms for solving RPCA, but these papers mainly focus on its convex relaxations. In this {paper}, we aim to review existing algorithms for both convex and nonconvex relaxations/variants of RPCA models and point out a few possible future directions.}


%
%

{
\begin{table*}[t]
  {\centering
  \caption{{Convergence rates of different algorithms.}}\label{tab:convergence-rate}
   \begin{tabular}{llll}\thickhline
    \textbf{Algorithm}                                      & \textbf{Problem}                 & $\epsilon$-\textbf{optimality measure} & \textbf{Convergence rate}         \\\thickhline
    \multicolumn{4}{c}{\textbf{Algorithms for Convex Models}} \\\hline
    PGM \cite{Lin-WCAMSAP-2009}                    & \eqref{rpca-nuclear-L1-noise-reform-unconstrained} & objective value error  & $\cO(1/k)$            \\\hline
    APGM \cite{Lin-WCAMSAP-2009}                   & \eqref{rpca-nuclear-L1-noise-reform-unconstrained} & objective value error  & $\cO(1/k^2)$        \\\hline
    IALM \cite{Lin-Chen-Wu-Ma-iadm-2009}           & \eqref{rpca-nuclear-L1} & --- & convergence, no rate given \\\hline
    ADMM \cite{Yuan-Yang-pjo-2013}                 & \eqref{rpca-nuclear-L1} & --- & convergence, no rate given \\\hline
    ALM \cite{Goldfarb-Ma-Scheinberg-2013}         & \eqref{rpca-nuclear-L1-smooth-g} & objective value error & $\cO(1/k)$ \\\hline
    FALM \cite{Goldfarb-Ma-Scheinberg-2013}        & \eqref{rpca-nuclear-L1-smooth-g} & objective value error & $\cO(1/k^2)$ \\\hline
    ASALM \cite{Tao-Yuan-siopt-2011}               & \eqref{rpca-nuclear-L1-noise-N-rewrite} & --- & convergence unclear, no rate given \\\hline
    VASALM \cite{Tao-Yuan-siopt-2011}               & \eqref{rpca-nuclear-L1-noise-N-rewrite} & --- & convergence, no rate given \\\hline
    PSPG \cite{Aybat-Ma-Goldfarb-2014}             & \eqref{rpca-nuclear-L1-noise}    & objective value error & $\cO(1/k)$ \\\hline
    ADMIP \cite{Aybat-Iyengar-2015}                & \eqref{rpca-nuclear-L1-noise}    & {objective value error} & {$\cO(1/k)$} \\\hline
    Quasi-Newton method {(fastRPCA)}~\cite{P1C2-Aravkin14_1P}       & \eqref{eq:aravkin_model}         & --- & {convergence,} no rate given \\\hline
    3-block ADMM \cite{Lin-Ma-Zhang-2015-free-gamma} & \eqref{lin-ma-zhang-spcp}      & --- & convergence, no rate given \\\hline
    Frank-Wolfe \cite{Mu-FW-2016}                  & \eqref{FW-RPCP}                  & objective value error & $\cO(1/k)$ \\\hline
        \multicolumn{4}{c}{\textbf{Algorithms for Nonconvex Models}} \\\hline
    GoDec \cite{GoDec-icml-2011} & \eqref{zhou-tao} & ---  & local convergence, no rate given \\\hline
    GreBsmo \cite{pmlr-v31-zhou13b} & \eqref{eq:biliteral-model} & --- & convergence unclear, no rate given \\\hline
    Alternating Minimization (R2PCP) \cite{hintermuller2015robust} & \eqref{eq:ams-model} & --- & local convergence, no rate given \\\hline
    Gradient Descent (GD) \cite{Yi-Park-Chen-Caramanis-gradient} & $\approx$ \eqref{eq:biliteral-model} & --- & linear convergence with proper initialization \\ &&& and incoherence assumption \\\hline
    Alternating Minimization \cite{Gu-Wang-Liu-AISTATS-2016} & \eqref{gu-wang-liu} & --- & local convergence with proper initialization \\ &&& and incoherence and RIP assumptions \\\hline
    Stochastic alg. \cite{Xu-OR-PCA-2013} & \eqref{OR-PCA} & --- & convergence if the iterates are always full rank matrices, \\ &&& no rate given \\\hline
    LMafit \cite{Shen-Wen-Zhang-Lmafit} & \eqref{Lmafit-problem} & --- & convergence if difference between two consecutive iterates \\ &&& tends to zero, no rate given \\\hline
    Conditional Gradient \cite{Jiang-Lin-Ma-Zhang-nonconvex-2016} & \eqref{nonconvex-rpca-cg-constrained} & perturbed KKT & $\cO(1/\sqrt{k})$ \\\hline
    ADMM \cite{Jiang-Lin-Ma-Zhang-nonconvex-2016} & \eqref{nonconvex-rpca-admm} & perturbed KKT & $\cO(1/\sqrt{k})$ \\\hline
    Proximal BCD \cite{Jiang-Lin-Ma-Zhang-nonconvex-2016} & \eqref{nonconvex-rpca-bcd} & perturbed KKT & $\cO(1/\sqrt{k})$ \\\thickhline
    \end{tabular}
    }

    \rule{0in}{1.2em}$^\dag$\scriptsize {Note: Some of these algorithms solve different problems and the $\epsilon$-optimality measures are also different, so the convergence rates are not directly comparable with each other. \cite{Yi-Park-Chen-Caramanis-gradient} has no explicit optimization formulation but the objective is similar to \eqref{eq:biliteral-model}. Moreover, global convergence is usually guaranteed for convex solvers, but only local convergence is usually guaranteed for nonconvex solvers, unless certain very strong assumptions are made.}
\end{table*}
}

\section{Algorithms for Convex Relaxations/Variants of RPCA}\label{sec:convex}

The
{earliest first-order methods for} solving the convex RPCP problem
\sa{are given in}~\cite{Lin-WCAMSAP-2009,Lin-Chen-Wu-Ma-iadm-2009}. In \cite{Lin-WCAMSAP-2009}, the authors 
proposed an accelerated proximal gradient method (APGM) \cite{Nesterov-07,Beck-Teboulle-2009,Tseng-2008} for solving \eqref{rpca-nuclear-L1-noise-reform-unconstrained} in which each iteration involves computing the proximal mappings of the nuclear norm $\|L\|_*$ and the $\ell_1$ norm $\|S\|_1$. In particular, the non-accelerated proximal gradient method (PGM) for solving \eqref{rpca-nuclear-L1-noise-reform-unconstrained} simply updates $L$ and $S$ as
\BE\label{pgm}
\BA{ll}
G^k     & := {\mu}(L^k+S^k-M) \\
L^{k+1} & := \argmin_L \ \|L\|_* + \frac{1}{2\tau}\|L-(L^k-\tau G^k)\|_F^2 \\
S^{k+1} & := \argmin_S \ \rho\|S\|_1 + \frac{1}{2\tau}\|S-(S^k-\tau G^k)\|_F^2.
\EA
\EE
Note that $G^k$ is the gradient of the quadratic penalty function in \eqref{rpca-nuclear-L1-noise-reform-unconstrained} and $\tau>0$ denotes a step size. The two subproblems in \eqref{pgm} both admit easy closed-form optimal solutions. Specifically, the solution of the $L$-subproblem corresponds to the proximal mapping of the nuclear norm, which is given by
\[L^{k+1} = \mathcal{S}_\tau(L^k-\tau G^k),\]
where the matrix shrinkage operation is defined as
\BE\label{mat-shrink}\mathcal{S}_\nu(Z) = U\diag((\sigma-\nu)_+)V^\top,\EE
where $Z = U\diag(\sigma)V^\top$ is the SVD of $Z$, and $z_+ := \max(0,z)$. The solution of the $S$-subproblem corresponds to the proximal mapping of the $\ell_1$ norm, which is given by
\[S^{k+1} = s_{\rho\tau}(S^k-\tau G^k),\]
where the vector shrinkage operation is defined as
\BE\label{vec-shrink}[s_\nu(Z)]_{ij} = \sign(Z_{ij})\circ\max\{0,|Z_{ij}|-\nu\},\EE
where $\sign(a)$ denotes the sign of $a$, and $\circ$ denotes the Hadamard product.


APGM incorporates Nesterov's acceleration technique and updates the variables as follows starting with $t_{-1}=t_0=1$.
\BE\label{apgm}
\BA{ll}
\bar{L}^k & := L^k + \frac{t_{k-1}-1}{t_k}(L^k-L^{k-1}) \\
\bar{S}^k & := S^k + \frac{t_{k-1}-1}{t_k}(S^k-S^{k-1}) \\[0.75mm]
\bar{G}^k & := {\mu}(\bar{L}^k+\bar{S}^k-M) \\
L^{k+1} & := \argmin_L \ \|L\|_* + \frac{1}{2\tau}\|L-(\bar{L}^k-\tau \bar{G}^k)\|_F^2 \\
S^{k+1} & := \argmin_S \ \rho\|S\|_1 + \frac{1}{2\tau}\|S-(\bar{S}^k-\tau \bar{G}^k)\|_F^2 \\
t_{k+1} & := (1+\sqrt{1+4t_k^2})/2.
\EA
\EE
Results in \cite{Nesterov-07,Beck-Teboulle-2009,Tseng-2008} show that the proximal gradient method \eqref{pgm} and the accelerated proximal gradient method \eqref{apgm} find an $\epsilon$-optimal solution to \eqref{rpca-nuclear-L1-noise-reform-unconstrained} in no more than $O(1/\epsilon)$ and $O(1/\sqrt{\epsilon})$ iterations, respectively.

\sa{When there is no noise, i.e., $\sigma=0$, the problem of interest is \eqref{rpca-nuclear-L1}. The drawback of the above approach for solving the unconstrained version in~\eqref{rpca-nuclear-L1-noise-reform-unconstrained} is that} \eqref{rpca-nuclear-L1-noise-reform-unconstrained} is equivalent to \eqref{rpca-nuclear-L1} only when ${\mu}\rightarrow +\infty$. Therefore, for any fixed ${\mu}>0$, there is always a residual term which does not go to zero. To remedy this, the same group of authors \cite{Lin-Chen-Wu-Ma-iadm-2009} considered the augmented Lagrangian method (ALM) for solving \eqref{rpca-nuclear-L1}. By associating a Lagrange multiplier $\Lambda$ to the linear equality constraint, the augmented Lagrangian function of \eqref{rpca-nuclear-L1} can be written as
\[\cL_\beta(L,S;\Lambda) := \|L\|_* + \rho\|S\|_1 - \langle \Lambda, L+S-M \rangle + \frac{\beta}{2}\|L+S - M\|_F^2,\]
where $\beta>0$ is a penalty parameter. A typical iteration of ALM {iterates the updates} as follows\sa{:}
\begin{subequations}
\label{rpca-nuclear-L1-ALM}
\begin{align}
(L^{k+1},S^{k+1}) & := \argmin_{L,S} \cL_\beta(L,S;\Lambda^k) \label{eq:ALM-subproblem}\\
\Lambda^{k+1} & := \Lambda^k - \beta(L^{k+1}+S^{k+1}-M).
\end{align}
\end{subequations}
Note that the first step in \eqref{rpca-nuclear-L1-ALM} requires to minimize the augmented Lagrangian function with respect to $L$ and $S$ simultaneously, which usually is \sa{computationally very expensive and almost as hard as} 
solving the original problem \sa{in}~{\eqref{rpca-nuclear-L1}}. In \cite{Lin-Chen-Wu-Ma-iadm-2009}, the authors proposed both exact and inexact versions of ALM, where the former one solves the subproblems (almost) exactly and the latter one solves the subproblems inexactly 
\sa{according to a particular subproblem termination criterion.} 
Both the exact ALM and inexact ALM (IALM) employ some iterative algorithm for minimizing the augmented Lagrangian function until certain \sa{overall} stopping criterion is met, which may require many iterations and thus time consuming. 
Around the same time when \cite{Lin-WCAMSAP-2009,Lin-Chen-Wu-Ma-iadm-2009} appeared, the alternating direction method of multipliers (ADMM) was revisited and found very successful in solving signal processing and image processing problems \cite{Combettes-2007,Goldstein-Osher-split-Bregman-2009,Yang-Zhang-Yin-ADMM-2010,Yang-Zhang-admm-2011}. It was then found that RPCP \sa{in} \eqref{rpca-nuclear-L1} can be nicely solved by ADMM due to its special separable structure \cite{Yuan-Yang-pjo-2013,Goldfarb-Ma-Scheinberg-2013}.
\sa{The ADMM iterations for solving \eqref{rpca-nuclear-L1} take the following form:} 
\BEA\label{rpca-nuclear-L1-ADMM}
\BA{ll}
L^{k+1} & := \argmin_{L} \cL_\beta(L,S^k;\Lambda^k) \\
S^{k+1} & := \argmin_{S} \cL_\beta(L^{k+1},S;\Lambda^k) \\
\Lambda^{k+1} & := \Lambda^k - \beta(L^{k+1}+S^{k+1}-M).
\EA
\EEA
Comparing to ALM \sa{in} \eqref{rpca-nuclear-L1-ALM}, it is noted that ADMM splits the subproblem in \sa{\eqref{eq:ALM-subproblem} into} two \sa{smaller} subproblems that correspond to computing proximal mappings of $\|L\|_*$ and $\|S\|_1$, respectively. The ADMM \eqref{rpca-nuclear-L1-ADMM} is known as two-block ADMM as there are two block variables $L$ and $S$ and hence two subproblems \sa{are solved} in each iteration of the algorithm. 
It is now widely known that the two-block ADMM is a special case of the so-called Douglas-Rachford operator splitting method \cite{Gabay-Mercier-1976,Gabay-83,Fortin-Glowinski-1983,Lions-Mercier-79} applied to the dual problem, and
the two-block ADMM for solving convex problems globally converges for any penalty parameter $\beta>0$ \cite{Eckstein-Bertsekas-1992} and converges with a sublinear rate $O(1/k)$ (see, e.g., \cite{He-Yuan-rate-ADM-2012,Monteiro-Svaiter-2010a,Lin-Ma-Zhang-convergence-2014}).

The alternating linearization method (ALM) proposed by Goldfarb, Ma and Scheinberg \cite{Goldfarb-Ma-Scheinberg-2013} is shown to be equivalent to a symmetric version of ADMM \eqref{rpca-nuclear-L1-ADMM} with either $\|L\|_*$ or $\|S\|_1$ \sa{replaced with some suitable smooth approximation}. \sa{For instance, given $\nu>0$, define $g_\nu:\reals^{m\times n}\rightarrow\reals$ such that
\begin{align}
\label{eq:smooth-g}
g_{\nu}(S)=\max_{Z\in\reals^{m\times n}}\{\fprod{S,Z}-\frac{\nu}{2}\norm{Z}_F^2:\ \norm{Z}_{\infty}\leq \rho\},
\end{align}
and let $g(S)=\rho \norm{S}_1$. Clearly, $g_\nu\rightarrow g$ uniformly as $\nu\searrow 0$. Moreover, $g_\nu$ is a differentiable convex function such that $\grad g_\nu$ is Lipschitz continuous. Indeed, given $S\in\reals^{m\times n}$, let $Z_\nu(S)$ be the maximizer for \eqref{eq:smooth-g}, which in closed form can be written as $Z_\nu(S)={\sign(S)\circ}\max\{\frac{1}{\nu}|S|,~\rho \mathbf{1}_{m\times n}\}$, and $\grad g_\nu(S)=Z_\nu(S)$ is Lipschitz {continuous} with constant $C_{g_{\nu}}=\frac{1}{\nu}$. Similarly, given $\mu>0$, define $f_\mu:\reals^{m\times n}\rightarrow\reals$ such that
\begin{align}
\label{eq:smooth-f}
f_{{\mu}}(L)=\max_{W\in\reals^{m\times n}}\{\fprod{L,W}-\frac{\mu}{2}\norm{W}_F^2:\ \norm{W} \leq 1\},
\end{align}
where $\norm{\cdot}$ denotes the spectral norm, and let $f(L)=\lambda \norm{L}_*$. Clearly, $f_\mu\rightarrow f$ uniformly as $\mu\searrow 0$. Moreover, $f_\mu$ is a differentiable convex function such that $\grad f_\mu$ is Lipschitz continuous. Indeed, given $L\in\reals^{m\times n}$, let $W_\mu(L)$ be {the maximizer for \eqref{eq:smooth-f},} which in closed form can be written as \magenta{$W_\mu(L)=U\diag\big(\max\big\{\frac{\sigma}{\mu}-\mathbf{1},~\mathbf{0}\big\}\big)V^\top$,} where $L=U \diag(\sigma)V^\top$ is the singular value decomposition of $L$ with $\sigma\in\reals^r_{++}$ denoting the vector of singular values; moreover, $\grad f_\mu(S)=W_\mu(L)$ is Lipschitz {continuous} with constant $C_{f_{\mu}}=\frac{1}{\mu}$.}

The alternating linearization method in \cite{Goldfarb-Ma-Scheinberg-2013} can be applied to solve the following problem, which is a smoothed version of \eqref{rpca-nuclear-L1}:
\BE\label{rpca-nuclear-L1-smooth-g} \min \ \|L\|_* + g_{\nu}(S), \ \st, \ \ L+S = M.\EE
Denote the augmented Lagrangian function of \eqref{rpca-nuclear-L1-smooth-g} as $\tilde{\cL}_\beta(L,S;\Lambda)$, the alternating linearization method in  \cite{Goldfarb-Ma-Scheinberg-2013} {iterates the updates} as follows:
\BE\label{alternating-linearization-method}
\BA{ll}
L^{k+1} & := \argmin_{L} \tilde{\cL}_\beta(L,S^k;\Lambda^k), \\
\Lambda^{k+\half} & := \Lambda^k - \beta(L^{k+1}+S^{k}-M), \\
S^{k+1} & := \argmin_{S} \tilde{\cL}_\beta(L^{k+1},S;\Lambda^{k+\half}), \\
\Lambda^{k+1} & := \Lambda^{k+\half} - \beta(L^{k+1}+S^{k+1}-M).
\EA
\EE
The authors in \cite{Goldfarb-Ma-Scheinberg-2013} proved that this method has a sublinear convergence rate $O(1/k)$. They also proposed an accelerated version of \eqref{alternating-linearization-method} (FALM) in \cite{Goldfarb-Ma-Scheinberg-2013} by adopting Nesterov's acceleration technique, and proved that the accelerated alternating linearization method has a better sublinear convergence rate $O(1/k^2)$.

Based on the success of two-block ADMM for solving \eqref{rpca-nuclear-L1}, it is then very natural to apply ADMM to solve SPCP in \eqref{rpca-nuclear-L1-noise}. {To do so, one has to introduce a new variable $N$, and rewrite \eqref{rpca-nuclear-L1-noise} equivalently as
\BE\label{rpca-nuclear-L1-noise-N-rewrite}
\BA{ll}
\min_{L,S,N\in\mR^{m\times n}} & \|L\|_* + \rho\|S\|_1 \\ \st, & L+S+N=M, \norm{N}_F\leq\sigma.
\EA
\EE
The ADMM for solving \eqref{rpca-nuclear-L1-noise-N-rewrite} {iterates} as follows with three block variables:
\BEA\label{spcp-ADMM}
\BA{ll}
L^{k+1} & := \argmin_{L} \cL_\beta(L,S^k,N^k;\Lambda^k) \\
S^{k+1} & := \argmin_{S} \cL_\beta(L^{k+1},S,N^k;\Lambda^k) \\
N^{k+1} & := \argmin_{N} \cL_\beta(L^{k+1},S^{k+1},N;\Lambda^k) \\
\Lambda^{k+1} & := \Lambda^k - \beta(L^{k+1}+S^{k+1}+N^{k+1}-M),
\EA
\EEA
where the augmented Lagrangian function for \eqref{rpca-nuclear-L1-noise-N-rewrite} is defined as}
\sa{
\[\BA{l}\cL_\beta(L,S,N;\Lambda) := \|L\|_* + \rho\|S\|_1 + \Ind(N\mid\|N\|_F\leq\sigma) \\ - \langle \Lambda, L+S+N-M \rangle + \frac{\beta}{2}\|L+S+N-M\|_F^2,\EA\]}%
where $\Ind(N\mid\cN)$ denotes the indicator function of the set $\{N\in\cN\}$, i.e., $\Ind(N\mid\cN)=0$ if $N\in\cN$ and $\Ind(N\mid\cN)=\infty$ otherwise. Note that the three subproblems in \eqref{spcp-ADMM} all have closed-form solutions. In particular, the $L$-subproblem corresponds to the proximal mapping of $\|L\|_*$, the $S$-subproblem corresponds to proximal mapping of $\|S\|_1$, and the $N$-subproblem corresponds to projection onto the set $\{N\mid \|N\|_F\leq\sigma\}$. Similar idea was used in \cite{Peng-rasl-2012} for robust image alignments. In practice, this three-block ADMM usually works very well. However, it was later discovered that the ADMM with more than two block variables is not necessarily convergent in general \cite{Chen-admm-failure-2013}. \sa{Note that although {\eqref{rpca-nuclear-L1-noise-N-rewrite}} contains three block variables, it can be viewed as a two-block problem, if we group $S$ and $N$ as one (larger) block variable. One of the earliest methods for solving SPCP in \eqref{rpca-nuclear-L1-noise} and \eqref{rpca-nuclear-L1-noise-N-rewrite} is a three-block ADMM algorithm, ASALM, proposed by Tao and Yuan~\cite{Tao-Yuan-siopt-2011}, and although it does not have any convergence guarantees, it works well in practice; and slightly changing the update rule in ASALM leads to VASALM, of which iterate sequence converges to an optimal solution; but this comes at the cost of degradation in practical convergence speed when compared to ASALM -- {indeed, VASALM \cite{Tao-Yuan-siopt-2011} can be seen as a linearized version of two-block ADMM to solve \eqref{rpca-nuclear-L1-noise} with a convergence guarantee without any convergence rate result.} To remedy the shortcoming associated with the theoretical convergence of three-block ADMM, several other alternatives based on two-block ADMM were proposed~\cite{aybat2011fast,Aybat-Iyengar-2015,Aybat-Ma-Goldfarb-2014}.}

\sa{Aybat, Goldfarb and Ma \cite{Aybat-Ma-Goldfarb-2014} proposed an accelerated proximal gradient method, \sa{PSPG,} for solving SPCP in~\eqref{rpca-nuclear-L1-noise}. 
First, \eqref{rpca-nuclear-L1-noise} is reformulated with a partially smooth objective. In particular, the nuclear norm is smoothed according to \eqref{eq:smooth-f}:
\begin{align}
\label{rpca-nuclear-L1-noise-rewrite-smooth}
\min_{L,S\in\mR^{m\times n}} f_\mu(L) + \rho\|S\|_{1}\ \st\ (L,S)\in\chi,\\
\chi:=\{(L,S)\mid\ \|L+S - M\|_F \leq \sigma\},
\end{align}
where $\mu>0$ is a given smoothing parameter. An accelerated proximal gradient method such as~\cite{Beck-Teboulle-2009,Nesterov-07}} can be applied to solve \eqref{rpca-nuclear-L1-noise-rewrite-smooth}, because it was shown in \cite{Aybat-Ma-Goldfarb-2014} that the following subproblem is easy to solve:
\begin{align}
\label{eq:complicated-prox}
\min_{L,S} \|S\|_1+\frac{1}{2\xi}\|L-\tilde{L}\|_F^2\ \st\ (L,S)\in\chi,
\end{align}
where $\xi>0$ denotes a step size of the proximal gradient step and $\tilde{L}$ denotes some known matrix. \sa{This operation requires one sorting which has $\cO(mn\log(mn))$ complexity.}

\sa{For any $\epsilon>0$, setting $\mu=\Omega(\epsilon)$, PSPG proposed in~\cite{Aybat-Ma-Goldfarb-2014} can compute an $\epsilon$-optimal solution to \eqref{rpca-nuclear-L1-noise} within $\cO(1/\epsilon)$ iterations, and its computational complexity per iteration is comparable to the work per iteration required by ASALM and VASALM, which is mainly determined by an SVD computation. On the other hand, it is also important to emphasize that PSPG iterate sequences do not converge to an optimal solution to the SPCP problem in \eqref{rpca-nuclear-L1-noise}. In particular, since within PSPG the smoothing parameter $\mu$ is fixed, depending on the approximation parameter $\epsilon$ for solving \eqref{rpca-nuclear-L1-noise-rewrite-smooth}, further iterations after reaching an $\epsilon$-optimal solution in $\cO(1/\epsilon)$ iterations do not necessarily improve the solution quality.}

In~\cite{aybat2011fast,Aybat-Iyengar-2015}, the \emph{variable penalty} ADMM algorithm ADMIP (Alternating Direction Method with Increasing Penalty) is proposed to solve the following equivalent formulation for \eqref{rpca-nuclear-L1-noise} using the variable splitting trick:
\BE\label{rpca-nuclear-L1-noise-equiv}
\min_{\hat{L},L,S\in\mR^{m\times n}} \|\hat{L}\|_* + \rho\|S\|_1 \ \st\ (L,S)\in\chi,\quad L=\hat{L}.
\EE
The augmented Lagrangian function of \eqref{rpca-nuclear-L1-noise-equiv} can be written as
\[\cL_\beta(\hat{L},L,S;\Lambda) := \|\hat{L}\|_* + \rho\|S\|_1 - \langle \Lambda, \hat{L}-L \rangle + \frac{\beta}{2}\|\hat{L}-L\|_F^2.\]
Given a nondecreasing penalty parameter sequence $\{\beta^k\}_{k\in\integers_+}$, ADMIP {updates the variables} as follows\sa{:}
\begin{subequations}
\label{rpca-nuclear-L1-ADMIP}
\begin{align}
\hat{L}^{k+1} & := \argmin_{\hat{L}} \cL_{\beta^k}(\hat{L},L^k,S^k;\Lambda^k) \label{eq:ADMIP-subproblem-1}\\
(L^{k+1},S^{k+1}) & := \argmin_{(L,S)\in\chi} \cL_{\beta^k}(\hat{L}^{k+1},L,S;\Lambda^k) \label{eq:ADMIP-subproblem-2}\\
\Lambda^{k+1} & := \Lambda^k - \beta^k(\hat{L}^{k+1}-L^{k+1}).
\end{align}
\end{subequations}
The step in~\eqref{eq:ADMIP-subproblem-1} requires computing a soft thresholding on the singular values of an $m\times n$ matrix and the step in~\eqref{eq:ADMIP-subproblem-2} requires an operation given in \eqref{eq:complicated-prox}.

Under mild conditions on the penalty parameter sequence, Aybat and Iyengar show that the primal-dual ADMIP iterate sequence converges to an optimal primal-dual solution to the SPCP problem in \eqref{rpca-nuclear-L1-noise-equiv} -- hence, $\{(L^k,S^k)\}_{k\in\integers_+}$ converges to an optimal solution to \eqref{rpca-nuclear-L1-noise}, and when constant penalty parameter is used as a special case, it can compute an $\epsilon$-optimal solution within $\cO(1/\epsilon)$ iterations, of which complexity is determined by an SVD. In particular, one needs the penalty parameter sequence $\{\beta^k\}_{k \in\integers_+}$ to be non-decreasing and to satisfy $\sum_{k}(\beta^k)^{-1}=+\infty$. The main advantages of adopting an increasing sequence of penalties are as follows:
  \begin{enumerate}
  \item The algorithm is robust in the sense that there is no need to search for problem data dependent $\beta^*$ that works well in practice. 
  \item The algorithm is likely to achieve primal feasibility faster.
  \item The complexity of initial (transient) iterations can be controlled through controlling $\{\beta^k\}$. 
    The main computational bottleneck in ADMIP is the SVD computation in \eqref{eq:ADMIP-subproblem-1}. Since the optimal $L^*$ is of low-rank, and $L_k\rightarrow L^*$, eventually the SVD computations are likely to be very efficient. However, since the initial iterates in the transient phase of the algorithm may have large rank, the complexity of the SVD in the initial iterations
    can be quite large. To compute the solution to the subproblem in~\eqref{eq:ADMIP-subproblem-1}, one does not need to compute singular values smaller than $1/\beta^k$; hence, initializing ADMIP with a \emph{small} $\beta^0>0$ will significantly decrease the complexity of initial iterations through employing \emph{partial} SVD computations, e.g., Lanczos-based methods such as PROPACK~\cite{larsen1998lanczos}.
  \end{enumerate}
In~\cite{Aybat-Iyengar-2015}, Aybat and Iyengar compared ADMIP against ASALM on both randomly generated synthetic problems and surveillance video foreground extraction problems. According to numerical results reported in~\cite{Aybat-Iyengar-2015}, on the synthetic problems ASALM requires about \emph{twice} as many iterations for convergence, while the total runtime for ASALM is considerably larger.

\sa{Aravkin et al.~\cite{P1C2-Aravkin14_1P} proposed solving
\begin{equation}
\label{eq:aravkin_model}
\min\psi(L,S)\quad \hbox{ s.t. }\quad \phi(\cA(L,S)-M) \leq \sigma,
\end{equation}
where $\cA:\reals^{m\times n}\times\reals^{m\times n}\rightarrow\reals^{m\times n}$ is a linear operator, $\phi:\reals^{m\times n}\rightarrow\reals$ is a smooth convex loss-function, and $\psi$ can be set to either one of the following functions:
\begin{align*}
\psi_{\rm sum}(L,S)&:=\norm{L}_* + \rho\norm{S}_1,\\
\psi_{\rm max}(L,S)&:=\max\{\norm{L}_*, \rho_{\max}\norm{S}_1\},
\end{align*}
$\rho,\rho_{\max}>0$ are some given function parameters. Note that setting $\psi=\psi_{\rm sum}$, $\rho(.)=\norm{.}^2_F$, and $\cA(L,S)=\pi_{\Omega}{L+S}$ in \eqref{eq:aravkin_model}, one obtains the SPCP problem in \eqref{rpca-nuclear-L1-noise-partial}. This approach offers advantages over the original SPCP formulation in terms of practical parameter selection. The authors make a case that although setting $\rho=\frac{1}{\sqrt{\max\{m,n\}}}$ in \eqref{rpca-nuclear-L1-noise} has theoretical justification as briefly discussed in the introduction section, many practical problems may violate the underlying assumptions in \eqref{incoherence-conditions}; in those cases one needs to tune $\rho$ via cross validation, and selecting $\rho_{\max}$ in $\psi_{\max}$ might be easier than selecting $\rho$ in $\psi_{\rm sum}$. Instead of solving \eqref{eq:aravkin_model} directly, a convex variational framework, accelerated with a ``quasi-Newton" method, is proposed. In particular, Newton's method is used to find a root of the value function:
\BE\label{eq:aravkin_value_function}
\BA{ll}
\upsilon(\tau):= & \min_{L,S\in\reals^{m\times n}}\phi(\cA(L,S)-A)-\sigma \\
                 & \st, \psi(L,S)\leq\tau,
\EA
\EE
i.e., given $\sigma>0$ compute $\tau^*=\tau(\sigma)$ such that $\upsilon(\tau^*)=0$. According to results in~\cite{P1C2-Aravkin13_1J}, if the constraint in \eqref{eq:aravkin_model} is tight at an optimal solution, then there exists $\tau^*=\tau(\sigma)$ such that $\upsilon(\tau^*)=0$ and the corresponding optimal solution to \eqref{eq:aravkin_value_function} is also optimal to \eqref{eq:aravkin_model}. Within Newton's method for root finding, to compute the next iterate $\tau^{k+1}$, one can compute the derivative of the value function at the current iterate $\tau^k$ as follows $\upsilon'(\tau_k)=-\psi^\circ(\cA^\top\grad\phi(\cA(L_k,S_k)-A))$, where $\psi^\circ$ denotes the polar gauge to $\psi$ and $(L_k,S_k)$ denotes the optimal solution to \eqref{eq:aravkin_value_function} at $\tau=\tau_k$ -- Aravkin et al. proposed a projected ``Quasi-Newton" method to solve \eqref{eq:aravkin_value_function}. According to numerical tests reported in \cite{P1C2-Aravkin13_1J}, QN-max, the quasi-Newton method running on \eqref{eq:aravkin_model} with $\psi=\psi_{\rm max}$ and $\phi(.)=\norm{.}_F^2$, is competitive with the state-of-the-art codes, ASALM~\cite{Tao-Yuan-siopt-2011}, PSPG~\cite{Aybat-Ma-Goldfarb-2014}, and ADMIP~\cite{Aybat-Iyengar-2015}.}

{
In a recent work \cite{Lin-Ma-Zhang-2015-free-gamma}, Lin, Ma and Zhang considered \sa{the penalty formulation of the SPCP problem, which is equivalent to solving \eqref{rpca-nuclear-L1-noise} for certain noise level $\sigma>0$:}
\BE\label{lin-ma-zhang-spcp}
\BA{ll}
\min & \|L\|_* + \rho\|S\|_1 + {\mu}\|\sa{N}\|_F^2 \\
\st  & L + S + \sa{N} = M,
\EA
\EE
where $\rho>0$ is the sparsity tradeoff parameter and $\mu>0$ is a suitable penalty parameter depending on the noise level $\sigma>0$. The authors showed that the following 3-block ADMM for solving \eqref{lin-ma-zhang-spcp} globally converges for any penalty parameter ${\beta}>0$.
{
\BE\label{lin-ma-zhang-spcp-admm}
\BA{ll}
L^{k+1} & := \argmin_L \ \cL_\beta(L,S^k,N^k;\Lambda^k) \\
S^{k+1} & := \argmin_S \ \cL_\beta(L^{k+1},S,N^k;\Lambda^k) \\
N^{k+1} & := \argmin_N \ \cL_\beta(L^{k+1},S^{k+1},N;\Lambda^k) \\
\Lambda^{k+1} & := \Lambda^k - \beta(L^{k+1} + S^{k+1} + N^{k+1} - M),
\EA
\EE}
where the augmented Lagrangian function is
\sa{
\[\BA{l}\cL_\beta(L,S,N;\Lambda) := \|L\|_* + \rho\|S\|_1 + \mu\|N\|_F^2 \\ - \langle \Lambda, L + S + N - M\rangle + \frac{{\beta}}{2}\|L + S + N - M\|_F^2.\EA\]}
Note that the three subproblems in \eqref{lin-ma-zhang-spcp-admm} are all easy to solve. Specifically, the $L$-subproblem corresponds to the proximal mapping of nuclear norm $\|L\|_*$, the $S$-subproblem corresponds to the proximal mapping of $\ell_1$ norm, and the $N$-subproblem admits a very easy analytical solution.
}

The Frank-Wolfe method (aka conditional gradient method) \cite{FW-1956} was revisited recently for solving large-scale machine learning problems \cite{Jaggi-icml-2013,thesisMJ}. RPCA is a representative example that is suitable for Frank-Wolfe method. Note that algorithms discussed above usually involve computing the proximal mapping of the nuclear norm, which is given by an SVD in \eqref{mat-shrink}. Computing full SVD for a large matrix in every iteration can be very time consuming. In contrast, the Frank-Wolfe method deals with nuclear norm in a much simpler manner, which only computes the largest singular value of a matrix in each iteration. The Frank-Wolfe method for solving RPCA was proposed by Mu \etal in \cite{Mu-FW-2016}. The authors in \cite{Mu-FW-2016} considered the penalized variant of RPCA \eqref{rpca-nuclear-L1-noise-reform-unconstrained}. However, this problem cannot be directly solved by the Frank-Wolfe method, because the Frank-Wolfe method requires a bounded constraint set. Therefore, the authors further reformulated \eqref{rpca-nuclear-L1-noise-reform-unconstrained} to the following problem for properly chosen constants $\lambda_L$, $\lambda_S$, $U_L$ and $U_S$:
\BE\label{FW-RPCP}
\BA{ll}
\min  & \half\|L+S-M\|_F^2 + \lambda_L t_L + \lambda_S t_S, \\
\st,                         & \|L\|_*\leq t_L \leq U_L,\ \|S\|_1\leq t_S \leq U_S.
\EA
\EE
The Frank-Wolfe method {iterates the updates} as follows:
\BE\label{FW-alg}
\BA{ll}
G^k & := L^k + S^k -M \\
\sa{(d_L^k,d_{t_L}^k)} & := \argmin_{\|L\|_*\leq{t_L}\leq U_L} \ \langle G^k,~L\rangle + \lambda_L{t_L}\\
\sa{(d_S^k,d_{t_S}^k)} & := \argmin_{\|S\|_1\leq{t_S}\leq U_S} \ \langle G^k,~S\rangle + \lambda_S{t_S} \\
\gamma^k &:= 2/(k+2)\\
L^{k+1} & := (1-\gamma^k)L^k + \gamma^k d_L^k  \\
t_L^{k+1} & := (1-\gamma^k)t_L^k + \gamma^k d_{t_L}^k \\
S^{k+1} & := (1-\gamma^k)S^k + \gamma^k d_S^k  \\
t_S^{k+1} & = (1-\gamma^k)t_S^k + \gamma^k d_{t_S}^k.
\EA
\EE
It was shown in \cite{Mu-FW-2016} that the two minimization subproblems in \eqref{FW-alg} are easy to solve. In particular, solving the subproblem for $(d_L,d_{t_L})$ requires only to compute the largest singular value and its corresponding singular vector of an $m\times n$ matrix. This is a big saving compared with computing the full SVD as required for computing the proximal mapping of the nuclear norm. As a result, Frank-Wolfe method has better per-iteration complexity than the proximal gradient method and ADMM algorithms discussed above, and thus may have better scalability for very large-scale problems. 
\sa{On the other hand, as pointed out in~\cite{Mu-FW-2016}, one clear disadvantage of Frank-Wolfe method on \eqref{FW-RPCP} is that at every iteration only one entry of the sparse component is updated. This leads to very slow convergence in practice. Hence, Mu et al.~\cite{Mu-FW-2016} proposed combining Frank-Wolfe iterations with an additional proximal gradient step in $S$-block. In particular, they proposed after Frank-Wolfe iterate $(L_{k+1},\tilde{S}_{k+1})$ is computed, an extra proximal gradient step is computed and $S$-block is updated again. Moreover, the authors also showed that this hybrid method obtained by combining Frank-Wolfe and proximal gradient steps enjoys a sublinear convergence rate $O(1/k)$ similar to Frank-Wolfe method given in~\eqref{FW-alg}.}

As a special case of RPCA, one can consider that all columns of the low-rank matrix $L$ are identical. That is, the given matrix $M$ is a superposition of a special rank-one matrix $L$ and a sparse matrix $S$. This special RPCA finds many interesting applications in practice such as video processing \cite{Li-Ng-Yuan-NLAA-2015,Li-Wang-Hu-Cai-2014} and bioinformatics \cite{Ma-Sparcoc-2015}. For instance, in the background extraction of surveillance video, if the background is static, then the low-rank matrix $L$ that corresponds to the background should have identical columns. As a result, the background and foreground can be separated by solving the following convex program:
\BE\label{nuclear-norm-free}
\BA{ll}
\min_{x,S} & \|S\|_1 \\
\st & [x,x,\cdots,x] \circ \mathbf{E} + S = M,
\EA
\EE
where $[x,x,\cdots,x]$ denotes the $m\times n$ matrix with all columns being $x$, $\mathbf{E}$ denotes the $m\times n$ {matrix} with all ones. Note that the optimal $x$ of \eqref{nuclear-norm-free} corresponds to the static background for all frames and $S$ corresponds to the moving foreground. The advantage of \eqref{nuclear-norm-free} is that it does not involve nuclear norm. As a result, SVD can be avoided when designing algorithms for solving it which makes the resulting algorithms very efficient. Yang, Pong and Chen \cite{Yang-Pong-Chen-2017} adopted the similar idea and designed variants of ADMM algorithm for solving a more general model where the sparsity function of $S$ is allowed to be a nonconvex function. Convergence of the proposed ADMM was proved under the assumption of KL property \cite{Kurdyka-1998,Lojasiewicz-1963} being satisfied. {We will discuss these topics in more details in the next section.}



\section{Algorithms for Nonconvex Relaxations/Variants of RPCA}

In this section, we discuss nonconvex relaxations and variants of RPCA given in~\eqref{rpca-rank-L0} and algorithms for solving them. \sa{Some researchers} aim to (approximately) solve RPCA in~\eqref{rpca-rank-L0} directly without convexifying the rank function \sa{and/or} the $\ell_0$ norm. In~\cite{GoDec-icml-2011}, Zhou and Tao considered 
\sa{a} variant of \eqref{rpca-rank-L0}:
\BE\label{zhou-tao}
\min \ \|L+S-M\|_F^2\ \st\ \rank(L)\leq \tau_r, \|S\|_0\leq \tau_s,
\EE
where $\tau_r$ and $\tau_s$ are given parameters to control the \sa{rank} of $L$ and sparsity of $S$. The authors proposed the GoDec algorithm which alternatingly minimizes the objective function \sa{in one variable while fixing the other, which is a special case of alternating projection method analyzed in~\cite{Lewis-Malick-alt-proj-manifold}.} 
\sa{In particular,} a naive version of GoDec algorithm {iterates} as follows:
\BE\label{GoDec}
\BA{ll}
L^{k+1} & := \argmin_L \|L+S^k-M\|_F^2\ \ \st\ \ \rank(L)\leq \tau_r, \\
S^{k+1} & := \argmin_S \|L^{k+1}+S-M\|_F^2\ \ \st\ \ \|S\|_0\leq \tau_s.
\EA
\EE
The two subproblems correspond to two projections. Although the projection for $S$ is easy, the projection for $L$ requires \sa{computing a partial SVD}, which may be time consuming when the matrix size is large.
The authors proposed to use \sa{a low-rank approximation based on  bilateral random projections} to approximate this projection \sa{operation} which can significantly speed up the computation. The authors showed that \sa{the iterate sequence} 
converges to a local minimum provided that the initial point is close to some point in the intersection of the two manifolds $\{L\mid\rank(L)\leq \tau_r\}$ and $\{S\mid \|S\|_0\leq \tau_s\}$. \sa{The convergence of GoDec 
follows from the results in~\cite{Lewis-Malick-alt-proj-manifold}.}

\sa{In~\cite{hintermuller2015robust}, Hinterm\"{u}ller and Wu considered a regularized version of \eqref{zhou-tao}:
\begin{align}
\label{eq:ams-model}
&\min_{L,S\in\reals^{m\times n}} \norm{L+S-M}_F^2+\frac{\rho}{2}\norm{L}_F^2\\
&\quad \hbox{ s.t. }\quad \rank(L)\leq \tau_{r},\quad \|S\|_0\leq\tau_{s},\nonumber
\end{align}
where $\tau_{r},\tau_{s}>0$ are given model parameters as in \eqref{zhou-tao}, and $0\leq\rho\ll 1$ is a given regularization parameter. An \emph{inexact alternating minimization method} (
\emph{R2PCP}) on matrix manifolds is proposed to solve \eqref{eq:ams-model}. The iterates $L_{k+1}$ and $S_{k+1}$ are computed as ``inexact" solutions to subproblems $\min_L\{\norm{L+S_k-M}_F^2+\rho\norm{L}_F^2:\ \rank(L)\leq\tau_{r}\}$ and $\min_S\{\norm{L_{k+1}+S-M}_F^2:\ \|S\|_0\leq\tau_{s}\}$, respectively. Provided that a limit point of the iterate sequence exists, under some further restrictive technical assumptions, it is shown that first-order necessary optimality conditions are satisfied.}

Note that \sa{the convex relaxation in~\eqref{rpca-nuclear-L1} involves the nuclear norm $\|\cdot\|_*$ in the objective}. Algorithms dealing with nuclear norm (like the ones discussed in \sa{Section~\ref{sec:convex}}) usually require to compute its proximal mapping, which then require an SVD. This can be very time consuming when the problem size is large, \sa{even when $\min\{m,~n\}$ is in the order of thousands.} This has motivated \sa{researchers} to consider nonconvex relaxations of RPCA that avoid SVD calculations. One way to achieve \sa{SVD-free methods} is to factorize the low rank matrix $L\in\mR^{m\times n}$ as \sa{a} product of \sa{two low-rank matrices,} 
i.e., factorize $L= UV^\top$, where $U\in\mR^{m\times r}$, $V\in\mR^{n\times r}$ and $r\ll \min\{m,n\}$ \sa{such that $r$ is an upper bound on $\rank(L^o)$.} This leads to many different nonconvex relaxations of RPCA.

{In~\cite{pmlr-v31-zhou13b}, Zhou and Tao considered a regularized version of \eqref{zhou-tao}:
\begin{align}
\label{eq:biliteral-model}
\min_{U,V,S} \rho_1\norm{S}_1 + \norm{UV^\top+S-M}_F^2\\
\quad \hbox{ s.t. }\quad \rank(U) = \rank(V) \leq \tau_{r},\nonumber
\end{align}
where $\rho_1>0$ and $\tau_{r}\in\integers_+$ such that $\tau_{r}\geq\rank(L^o)$. The authors propose a three-block alternating minimization algorithm, GreBsmo, for solving \eqref{eq:biliteral-model}. The proposed algorithm lacks theoretical convergence guarantees; but, on the other hand, according to numerical results reported in~\cite{pmlr-v31-zhou13b}, GreBsmo performs considerably better than both GoDec~\cite{GoDec-icml-2011} and inexact ALM method~\cite{Lin-Chen-Wu-Ma-iadm-2009} (around 30-100 times faster than both) when applied to foreground extraction problems.}

A nonconvex model of RPCA, {similar to one in~\cite{pmlr-v31-zhou13b},} is studied in \cite{Yi-Park-Chen-Caramanis-gradient} and a gradient descent (GD) method with proper initialization is proposed to solve it. The algorithm proposed in \cite{Yi-Park-Chen-Caramanis-gradient} has two phases and in both phases the objective is to reduce the function $Q(U,V,S):=\|UV^\top+S-M\|_F^2$. In the first phase, a sorting-based sparse estimator is used to generate a rough \sa{initial} estimate $S_0$ \sa{to the unknown sparse target matrix $S^o$}, and then $U_0$ and $V_0$ are generated via an SVD of $M-S_0$ such that $U_0V_0^\top$ \sa{forms a rough initial} estimate \sa{to the unknown low-rank target matrix $L^o$}. In the second phase, the algorithm alternatingly performs two operations: taking gradient steps for $U$ and $V$, and \sa{computing} a sparse estimator to adjust $S$. The sparse estimator is to guarantee that \sa{the fraction of nonzero entries in each column and row of $S$ is bounded above so that the nonzero entries are spread out in $S$}. The authors showed that the proposed two-phase algorithm \sa{recovers the target decomposition} and \sa{linear convergence is achieved} with proper initialization and step size, under the incoherence assumptions similar to \eqref{incoherence-conditions}.
For more detailed description of the assumptions and the results, see \cite{Yi-Park-Chen-Caramanis-gradient}.

In \cite{Gu-Wang-Liu-AISTATS-2016}, \sa{assuming that the data matrix $M$ is observed indirectly through compressive measurements}, Gu, Wang and Liu considered the following variant of RPCA:
\BE\label{gu-wang-liu}
\min \ H(U,V,S):=\|A(UV^\top+S)-M\|_F^2\ \ \st\ \ \|S\|_0\leq \tau_s,
\EE
where $A$ is a sensing matrix. The alternating minimization algorithm proposed in \cite{Gu-Wang-Liu-AISTATS-2016} iterates as follows:
\BEA\label{gu-wang-liu-alt-min}
\BA{ll}
U^{k+1} := \argmin_U \ H(U,V^k,S^k) \\
V^{k+1} := \argmin_V \ H(U^{k+1},V,S^k) \\
S^{k+1} := \argmin_S \ H(U^{k+1},V^{k+1},S)\ \ \st\ \ \|S\|_0\leq \tau_s.
\EA
\EEA
It is noted that the $U$ and $V$ subproblems in \eqref{gu-wang-liu-alt-min} correspond to solving linear systems and the $S$-subproblem admits an easily computable closed-form solution. The authors showed that \sa{under incoherence assumption on $L^o$ and $A$ satisfying restricted isometry property~(RIP),} \eqref{gu-wang-liu-alt-min} converges globally. However, note that $A=I$ does not satisfy the RIP condition, and therefore the convergence is not guaranteed for RPCA problem \sa{in \eqref{rpca-rank-L0}} for which $A=I$.
Similar idea was also \sa{investigated} in~\cite{sparcs-2011}, assuming the RIP condition on the sensing matrix $A$, and thus does not \sa{apply to the} RPCA problem either.

In \cite{Xu-OR-PCA-2013}, the authors considered \sa{the scenario such that the columns of the data matrix $M$ are observed} in an online fashion. This is suitable for many real applications, \sa{e.g., 
in surveillance video background separation.}
To handle this problem, the authors proposed a stochastic algorithm, which solves a nonconvex variant of RPCA:
\BE\label{OR-PCA}
\min_{U,V,S} \half\|UV^\top+S-M\|_F^2 + \frac{\rho_1}{2}(\|U\|_F^2 + \|V\|_F^2) + \rho_2\|S\|_1,
\EE
where $\rho_1$ and $\rho_2$ are \sa{some weight} parameters. \sa{The formulation in~\eqref{OR-PCA} exploits the representation of 
the nuclear norm established in~\cite{recht2010guaranteed}. In particular, for any given $L\in\reals^{m\times n}$ such that $\rank(L)\leq r$, $\norm{L}_*$ can be computed as follows:
\BE\label{equiv-nuc-nonconvex}\|L\|_*:=\inf_{U\in\mR^{m\times r},V\in\mR^{n\times r}}\big\{\half\|U\|_F^2+\half\|V\|_F^2: UV^\top=L\big\}.\EE}
Form \eqref{equiv-nuc-nonconvex} we know that \eqref{rpca-nuclear-L1-noise-reform-unconstrained} is equivalent to
\BE\label{OR-PCA-constrained}
\min_{U,V,S} \frac{\mu}{2}\|UV^\top+S-M\|_F^2 + \half(\|U\|_F^2 + \|V\|_F^2) + \rho\|S\|_1.
\EE
As a result, \eqref{OR-PCA} is a nonconvex reformulation of the penalized variant of RPCA in~\eqref{rpca-nuclear-L1-noise-reform-unconstrained}. For given matrix $M=[M_1,\ldots,M_n]\in\mR^{m\times n}$, solving \eqref{OR-PCA} is the same as the following empirical risk minimization problem:
\BE\label{OR-PCA-ERM}
\min_{U} \frac{1}{n}\sum_{i=1}^n \ell(M_i,U) + \rho_1\|U\|_F^2,
\EE
where $\ell(M_i,U)$ is defined as
\BE\label{ERM-component-func}\BA{ll}\ell(M_i,U):= & \min_{V_i\in\mR^r,S_i\in\mR^m} \half\|UV_i+S_i-M_i\|_2^2 \\ & + \frac{\rho_1}{2}\|V_i\|_2^2 + \rho_2\|S_i\|_1.\EA\EE
The empirical risk minimization \eqref{OR-PCA-ERM} favors stochastic gradient descent algorithm. Of course every time to compute the gradient of $\ell(M_i,U)$, another minimization problem in \eqref{ERM-component-func} needs to be solved. Therefore, the algorithm proposed in \cite{Xu-OR-PCA-2013} is an alternating minimization method with subproblem for $U$ being solved using stochastic gradient descent. The authors showed that the proposed method converges to the correct low-dimensional subspace asymptotically under certain assumptions.

The following nonconvex variant of RPCA was proposed by Shen, Wen and Zhang in \cite{Shen-Wen-Zhang-Lmafit}:
\BE\label{Lmafit-problem}
\min_{U,V} \ \|UV^\top -M\|_1.
\EE
This simple reformulation can be viewed as a nonconvex reformulation of \eqref{rpca-nuclear-L1} but without any regularization terms on $U$ and $V$. In particular, \eqref{Lmafit-problem} can be rewritten as
\sa{
\BE\label{Lmafit-problem-rewrite}
\min \ \|S\|_1\ \ \st\ \ S + UV^\top = M.
\EE}
The authors in \cite{Shen-Wen-Zhang-Lmafit} proposed an ADMM algorithm (named LMafit) for solving \eqref{Lmafit-problem-rewrite}.
By associating a Lagrange multiplier $\Lambda$ to the constraint, the augmented Lagrangian function for \eqref{Lmafit-problem-rewrite} can be written as
\begin{align*}
\cL_\beta(U,V,S;\Lambda):= & \|S\|_1-\langle \Lambda,UV^\top + S-M \rangle \\ & +\frac{\beta}{2}\|UV^\top + S-M\|_F^2,
\end{align*}
where $\beta>0$ is a penalty parameter.
The nonconvex ADMM for solving \eqref{Lmafit-problem-rewrite} {iterates the updates} as follows:
\BEA\label{Lmafit-admm}
\BA{ll}
U^{k+1} & := \argmin_U \ \cL_\beta(U,V^k,S^k;\Lambda^k) \\
V^{k+1} & := \argmin_V \ \cL_\beta(U^{k+1},V,S^k;\Lambda^k) \\
S^{k+1} & := \argmin_S \ \cL_\beta(U^{k+1},V^{k+1},S;\Lambda^k) \\
\Lambda^{k+1} & := \Lambda^k - \beta(U^{k+1}{V^{k+1}}^\top + S^{k+1}-M).
\EA
\EEA
Note that \sa{all} three subproblems in \eqref{Lmafit-admm} are easy to solve. In particular, the $U$ and $V$-subproblems correspond to solving linear systems, and the \sa{$S$-subproblem} corresponds to the soft-shrinkage operation of the $\ell_1$ norm \eqref{vec-shrink}. However, this nonconvex ADMM lacks convergence guarantees.

In \cite{Jiang-Lin-Ma-Zhang-nonconvex-2016}, the authors studied some variants of \sa{the conditional gradient method and ADMM for} solving nonconvex and nonsmooth optimization problems. \sa{Consider a general nonconvex optimization problem:}
\BE\label{cg-general-problem}\min_x \ f(x) + r(x)\ \ \st\ \ x\in\cX,\EE
where $f$ is smooth and possibly nonconvex, $r(x)$ is convex and nonsmooth, $\cX$ is a convex compact set. Moreover, $f$ satisfies the so-called H\"{o}lder condition:
\[f(y) \leq f(x) + \nabla f(x)^\top(y-x) + \frac{\sa{\gamma}}{2}\|y-x\|_p^p, \forall x,y\in\cX,\]
where $p>1$ and $\sa{\gamma}>0$. The definition of $\epsilon$-stationary solution given in \cite{Jiang-Lin-Ma-Zhang-nonconvex-2016} is as follows.
\begin{defn}\label{def:eps-unconstrained}
$x\in\cX$ is called an $\epsilon$-stationary solution ($\epsilon\geq 0$) for \eqref{cg-general-problem} if the following holds:
\[\Psi_\cX(x):=\inf_{y\in\cX}\{\nabla f(x)^\top(y-x) + r(y)-r(x)\}\geq-\epsilon.\]
\end{defn}
The authors commented that this definition is stronger than the one used by Ghadimi \etal in \cite{glz13}.

Now we briefly discuss how to apply the algorithms analyzed in~\cite{Jiang-Lin-Ma-Zhang-nonconvex-2016} to solve nonconvex RPCA variants.
\sa{Consider the nonconvex RPCA variant given in \eqref{OR-PCA}, 
the generalized conditional gradient method proposed in \cite{Jiang-Lin-Ma-Zhang-nonconvex-2016} can be customized} to solve~\eqref{OR-PCA}. 
\sa{Since} the generalized conditional gradient method requires a compact constraint set, 
\sa{one can equivalently reformulate \eqref{OR-PCA} in the following form:}
\BE\label{nonconvex-rpca-cg-constrained}
\BA{ll}
\min_{U,V,S} & f(U,V,S) + \rho_2\|S\|_1, \\[1mm]
\st         & \sa{\max\{\|U\|_F,~\|V\|_F\}} \leq \tfrac{1}{\sqrt{\rho_1}}\|M\|_F, \\[1mm]
             & \|S\|_1 \leq \tfrac{1}{2\rho_2}\|M\|_F^2,
\EA
\EE
where $f(U,V,S):=\half\|UV^\top+S-M\|_F^2 + \frac{\rho_1}{2}(\|U\|_F^2 + \|V\|_F^2)$ denotes the smooth part of the objective function. It is easy to see that $\nabla f$ is Lipschitz continuous \sa{-- let $\gamma>0$ 
denote the Lipschitz constant.}

\sa{At the $k$-th iteration of the generalized conditional gradient method~\cite{Jiang-Lin-Ma-Zhang-nonconvex-2016}, implemented on \eqref{nonconvex-rpca-cg-constrained},
one needs to solve the following subproblem:} 
\BE\label{nonconvex-rpca-cg-constrained-sub}
\BA{ll}
\min_{U,V,S} & \langle \nabla_U f(U^k,V^k,S^k),U \rangle + \langle \nabla_V f(U^k,V^k,S^k),V \rangle \\
             & +\langle \nabla_S f(U^k,V^k,S^k),S \rangle  + \rho_2\|S\|_1, \\[1mm]
\st          & \max\{\|U\|_F,~\|V\|_F\} \leq \tfrac{1}{\sqrt{\rho_1}}\|M\|_F, \\[1mm]
             & \|S\|_1 \leq \tfrac{1}{2\rho_2}\|M\|_F^2.
\EA
\EE
Suppose $(\bar{U}^k,\bar{V}^k,\bar{S}^k)$ denotes the solution of \eqref{nonconvex-rpca-cg-constrained-sub}, a typical iteration of the generalized conditional gradient method {is given} as follows:
\BE\label{nonconvex-rpca-cg-alg}\left\{
\BA{ll}
\mbox{Solve } & \eqref{nonconvex-rpca-cg-constrained-sub} \mbox{ to obtain } (\bar{U}^k,\bar{V}^k,\bar{S}^k) \\
\Delta U^k & := \bar{U}^k - U^k \\
\Delta V^k & := \bar{V}^k - V^k \\
\Delta S^k & := \bar{S}^k - S^k \\
\Delta^k  &:= \sa{[\Delta U^k,\Delta V^k,\Delta S^k]} \\
\alpha_k & := \argmin_{\alpha\in[0,1]} \alpha (\langle \nabla f(U^k,V^k,S^k),\Delta^k \rangle \\ & + \frac{\alpha^2\sa{\gamma}}{2}\|\Delta^k\|_F^2 + (1-\alpha)\rho_2\|[U^k,V^k,S^k]\|_1 \\ & + \alpha\rho_2\|[\bar{U}^k,\bar{V}^k,\bar{S}^k]\|_1 \\
U^{k+1} & := (1-\alpha_k) U^k + \alpha_k \bar{U}^k \\
V^{k+1} & := (1-\alpha_k) V^k + \alpha_k \bar{V}^k \\
S^{k+1} & := (1-\alpha_k) S^k + \alpha_k \bar{S}^k.
\EA\right.
\EE
Note that the generalized conditional gradient method in \cite{Jiang-Lin-Ma-Zhang-nonconvex-2016} involves a line search step for computing $\alpha_k$ as shown in~\eqref{nonconvex-rpca-cg-alg}, \sa{which can be efficiently computed.} 
It is shown in \cite{Jiang-Lin-Ma-Zhang-nonconvex-2016} that the generalized conditional gradient method in~\eqref{nonconvex-rpca-cg-alg} \sa{can compute} an $\epsilon$-stationary solution of \eqref{nonconvex-rpca-cg-constrained} in $O(1/\epsilon^2)$ iterations.

Jiang~\etal \cite{Jiang-Lin-Ma-Zhang-nonconvex-2016} also proposed \sa{some ADMM variants that can solve various nonconvex RPCA formulations, and the authors provided a convergence rate analysis to compute an $\epsilon$-stationary solution -- the definition of $\epsilon$-stationarity employed to analyze the ADMM algorithm is given in Definition~\ref{def:eps-constrained}.} We here discuss the ADMM-g algorithm in \cite{Jiang-Lin-Ma-Zhang-nonconvex-2016} which can solve the following RPCA variant:
\BE\label{nonconvex-rpca-admm}
\BA{ll}
\min & \half\|L-UV^\top\|_F^2 + \frac{\rho_1}{2}(\|U\|_F^2 + \|V\|_F^2) \\[1mm] & + \rho_2\|S\|_1 + \rho_3\|N\|_F^2 \\[1mm]
\st  & L + S + N = M.
\EA
\EE
This is a nonsmooth and nonconvex problem with five block variables $L$, $S$, $N$, $U$ and $V$, and it can be viewed as a variant of \eqref{OR-PCA} with linear constraints. Treating $N$ as the last block variable, a typical iteration of ADMM-g for solving \eqref{nonconvex-rpca-admm} {iterates} as follows:
\BE\label{ADMM-g}
\BA{ll}
L^{k+1} &:= \argmin_L \ \tilde{\cL}_{L^k}(L,U^{k},V^{k},S^{k},N^{k};\Lambda^k) \\
U^{k+1} &:= \argmin_U \ \tilde{\cL}_{U^k}(L^{k+1},U,V^{k},S^{k},N^{k};\Lambda^k) \\
V^{k+1} &:= \argmin_V \ \tilde{\cL}_{V^k}(L^{k+1},U^{k+1},V,S^{k},N^{k};\Lambda^k) \\
S^{k+1} &:= \argmin_S \ \tilde{\cL}_{S^k}(L^{k+1},U^{k+1},V^{k+1},S,N^{k};\Lambda^k) \\
N^{k+1} &:= N^k - \sa{\eta}\nabla_N\cL(L^{k+1},U^{k+1},V^{k+1},S^{k+1},N;\Lambda^k) \\
\Lambda^{k+1} &:= \Lambda^k - \beta(L^{k+1}+S^{k+1}+N^{k+1}-M),
\EA
\EE
where $\sa{\eta}>0$ is a step size, the augmented Lagrangian function $\cL$ is defined as
\begin{eqnarray*}
\lefteqn{\cL(L,U,V,S,N;\Lambda)  := }\\
& & \half\|L-UV^\top\|_F^2 + \frac{\rho_1}{2}(\|U\|_F^2 + \|V\|_F^2) + \rho_2\|S\|_1 \\
& & + \rho_3\|N\|_F^2  -\langle\Lambda,L + S + N - M \rangle \\
& & + \frac{\beta}{2}\|L + S + N - M\|_F^2,
\end{eqnarray*}
and $\tilde{\cL}$ denotes $\cL$ plus a proximal term. For example, $\tilde{\cL}_{L^k}$ is defined as
\begin{eqnarray*}
\lefteqn{\tilde{\cL}_{L^k}(L,U^{k},V^{k},S^{k},N^{k};\Lambda^k) :=}\\
& & \cL(L,U^{k},V^{k},S^{k},N^{k};\Lambda^k) + \half\|L-L^k\|_H^2,
\end{eqnarray*}
where $H$ denotes a pre-specified positive definite matrix which needs to satisfy certain conditions to guarantee the convergence of \sa{the method as stated in~\eqref{ADMM-g}.}
It is noted that the last block variable $N$ is treated specially. It is not updated by minimizing the augmented Lagrangian function, but by taking a gradient step on it.
The results in \cite{Jiang-Lin-Ma-Zhang-nonconvex-2016} indicate that ADMM-g \eqref{ADMM-g} finds an $\epsilon$-stationary solution for \eqref{nonconvex-rpca-admm} in no more than $O(1/\epsilon^2)$ iterations. Since \eqref{nonconvex-rpca-admm} is a constrained problem, the definition of its $\epsilon$-stationary solution is different from the one in Definition \ref{def:eps-unconstrained}. Here we briefly discuss how it is defined for constrained problems in \cite{Jiang-Lin-Ma-Zhang-nonconvex-2016}. We consider the following constrained nonsmooth and nonconvex problem:
\BE\label{admm-general-problem}\BA{ll}
\min & f(x_1,\ldots,x_p) + \sum_{i=1}^{p-1}r_i(x_i) \\
\st  & \sum_{i=1}^{p-1} A_ix_i + x_p = b, \\
     & x_i\in\cX_i, i=1,\ldots,p-1,
\EA\EE
where $x_i\in\mR^{n_i}$, $f$ is differentiable and possibly nonconvex, each $r_i$ is possibly nonsmooth and nonconvex, and each $\cX_i$ is a convex set. Note again that the last block variable $x_p$ is treated differently, which is needed in the analysis of convergence rate. The $\epsilon$-stationary solution to \eqref{admm-general-problem} is defined as follows.
\begin{defn}\label{def:eps-constrained}
$(x_1^*,\ldots,x_p^*)\in\cX_1\times\cdots\times\cX_{p-1}\times\mR^{n_p}$ is called an $\epsilon$-stationary solution to \eqref{admm-general-problem}, if there exists $\lambda^*$ such that the following holds for any $(x_1,\ldots,x_p)\in\cX_1\times\cdots\times\cX_{p-1}\times\mR^{n_p}$:
\[\BA{r}
(x_i-{x}^*_i)^\top({g}^*_i + \nabla_i f(x_1^*,\cdots, x^*_p)  - A_i^\top {\lambda}^*) \geq -\epsilon, \\
 i=1,\ldots,p-1, \\
\| \nabla_p f(x_1^*,\ldots, x_{p-1}^*, x_p^*) - \lambda^*\| \leq  \epsilon, \\
\| \sum_{i=1}^{p-1} A_i x_i^* + x_p^* - b\| \leq \epsilon,
\EA\]
where ${g}^*_i$ is a general subgradient of $r_i$ at point $x_i^*$. This set of inequalities can be viewed as a perturbed KKT system.
\end{defn}

It is also interesting to note that \eqref{nonconvex-rpca-admm} is equivalent to the following unconstrained problem, and thus can be solved by block coordinate descent method (BCD).
\BE\label{nonconvex-rpca-bcd}
\BA{ll}
\min & \half\|M-S-N-UV^\top\|_F^2 + \frac{\rho_1}{2}(\|U\|_F^2 + \|V\|_F^2) \\ & + \rho_2\|S\|_1 + \rho_3\|N\|_F^2.
\EA
\EE
Most existing BCD type algorithms for solving nonconvex problems lack convergence rate analysis. In \cite{Jiang-Lin-Ma-Zhang-nonconvex-2016}, the authors proposed a proximal BCD method that can find an $\epsilon$-stationary solution in $O(1/\epsilon^2)$ iterations. Denoting the objective function in \eqref{nonconvex-rpca-bcd} as $F(U,V,S,N)$, the proximal BCD given in \cite{Jiang-Lin-Ma-Zhang-nonconvex-2016} for solving \eqref{nonconvex-rpca-bcd} {updates the variables} as follows:
\BE\label{prox-bcd}\BA{l}
U^{k+1} := \argmin_U F(U,V^k,S^k,N^k)+\half\|U-U^k\|_H^2 \\
V^{k+1} := \argmin_V F(U^{k+1},V,S^k,N^k)+\half\|V-V^k\|_H^2 \\
S^{k+1} := \argmin_S F(U^{k+1},V^{k+1},S,N^k)+\half\|S-S^k\|_H^2 \\
N^{k+1} := \argmin_N F(U^{k+1},V^{k+1},S^{k+1},N)+\half\|N-N^k\|_H^2,
\EA\EE
where $H$ denotes a pre-specified positive definite matrix. 

\section{Preliminary Numerical Experiments}

\begin{figure*}[t]
  \caption{Comparison of different algorithms. The first row shows the comparison results for $n=500$, $SNR=50dB$, and $(c_r,c_p) = (0.05,0.05)$, $(0.05,0.1)$, $(0.1,0.05)$, $(0.1,0.1)$, respectively from left to right. The second row corresponds to $n=500$, $SNR=100dB$, the third row corresponds to $n=1500$, $SNR=50dB$, and the fourth row corresponds to $n=1500$, $SNR=100dB$, all with $(c_r,c_p)$ set the same as the first row from left to right. We remark again that these comparison results only reflect the performance of the algorithms under the current setting of parameters, stopping criteria, initial points, input data and so on. For other settings and data, the performance can be very different. See Remark \ref{rem:comparison}.}\label{fig:synthetic}
  \centering
  \includegraphics[width=0.24\textwidth]{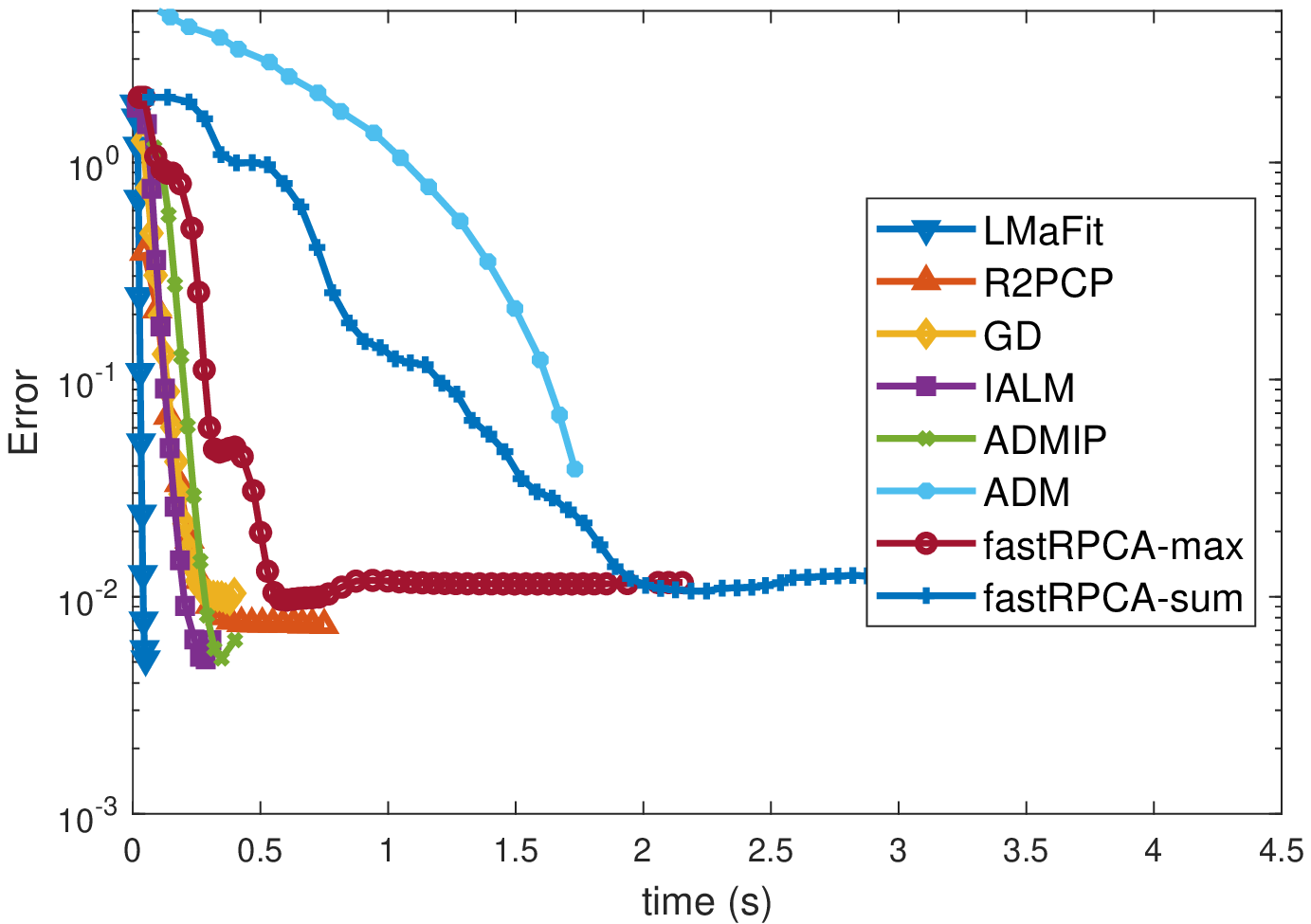}
  \includegraphics[width=0.24\textwidth]{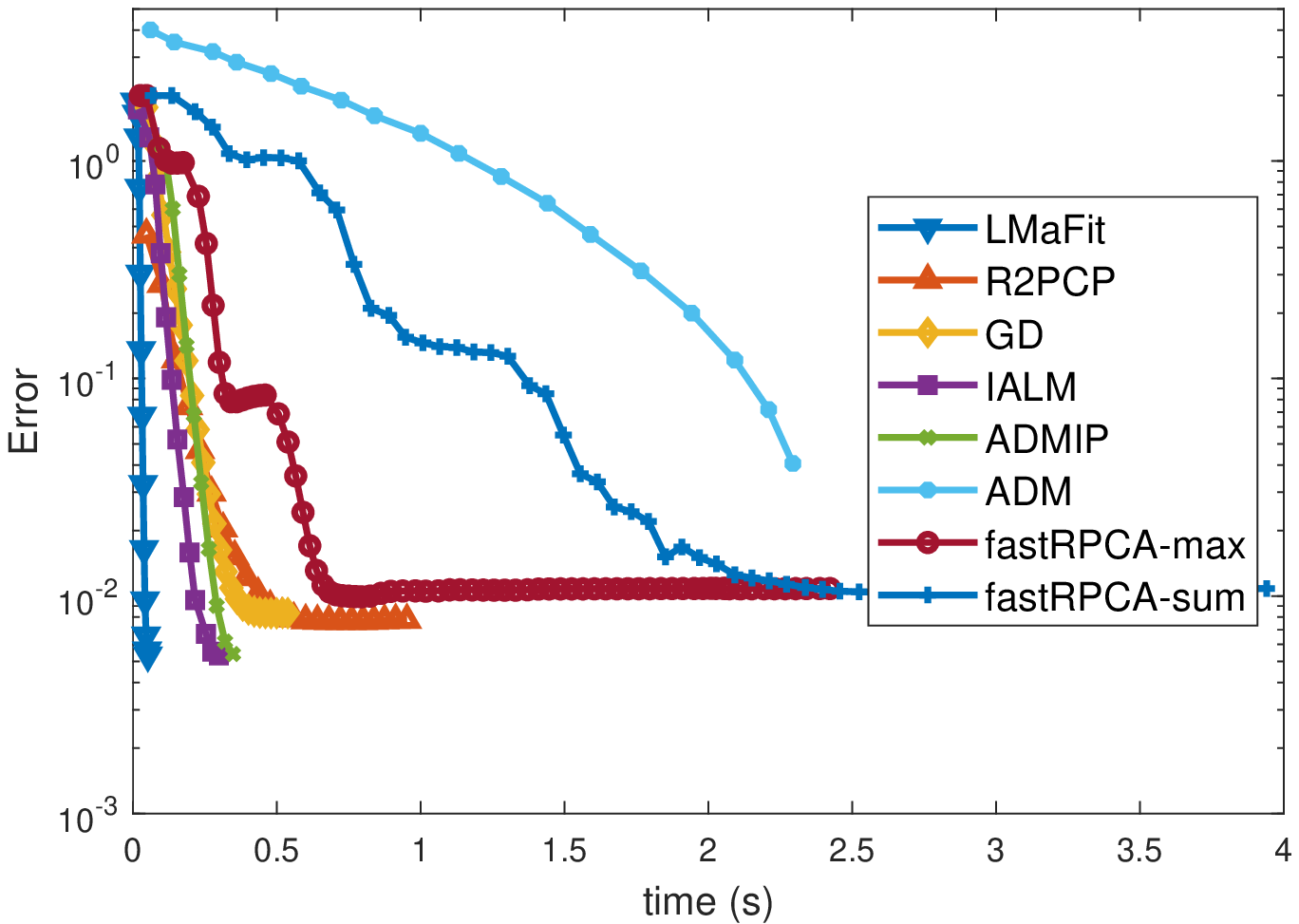}
  \includegraphics[width=0.24\textwidth]{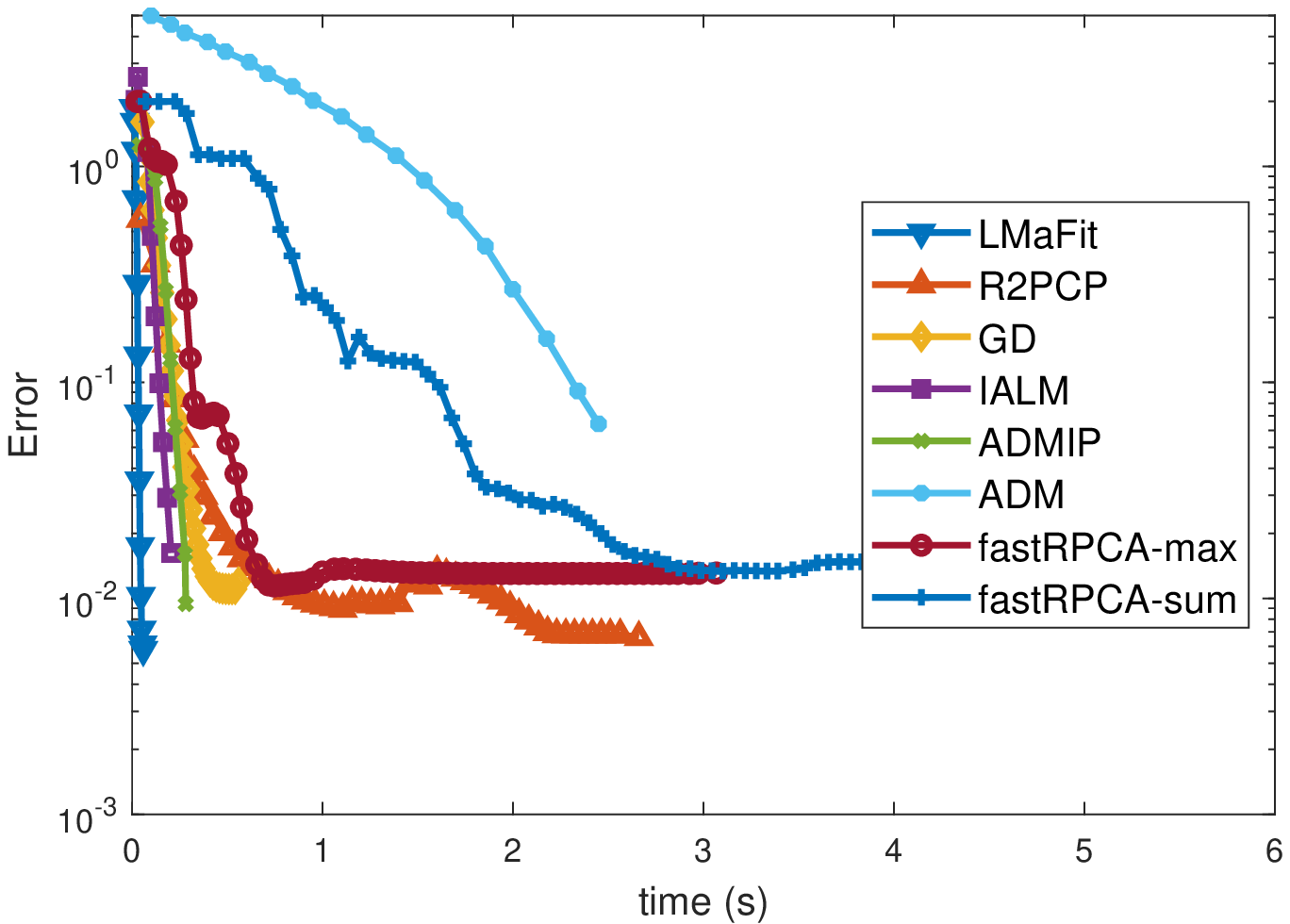}
  \includegraphics[width=0.24\textwidth]{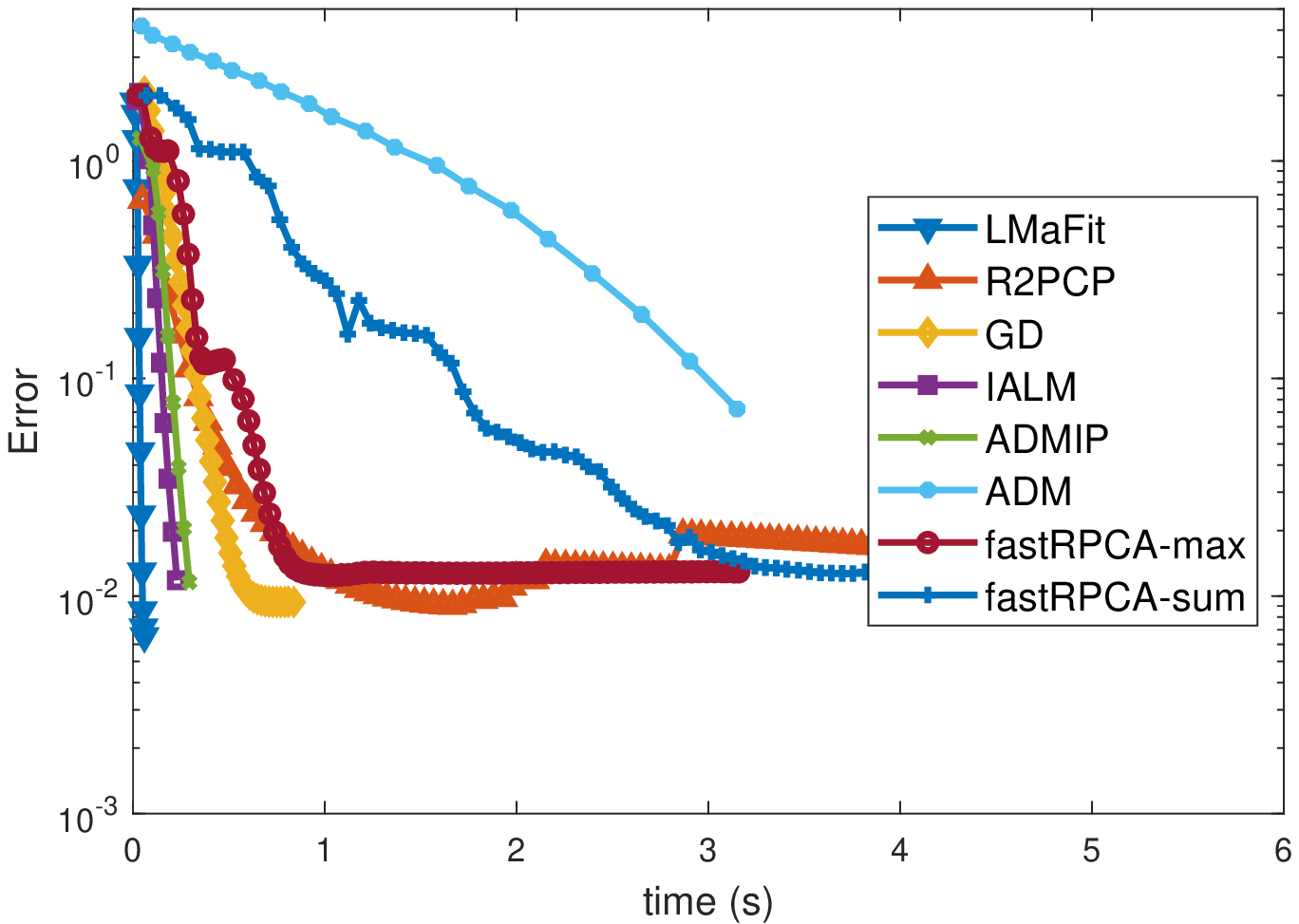}
  \includegraphics[width=0.24\textwidth]{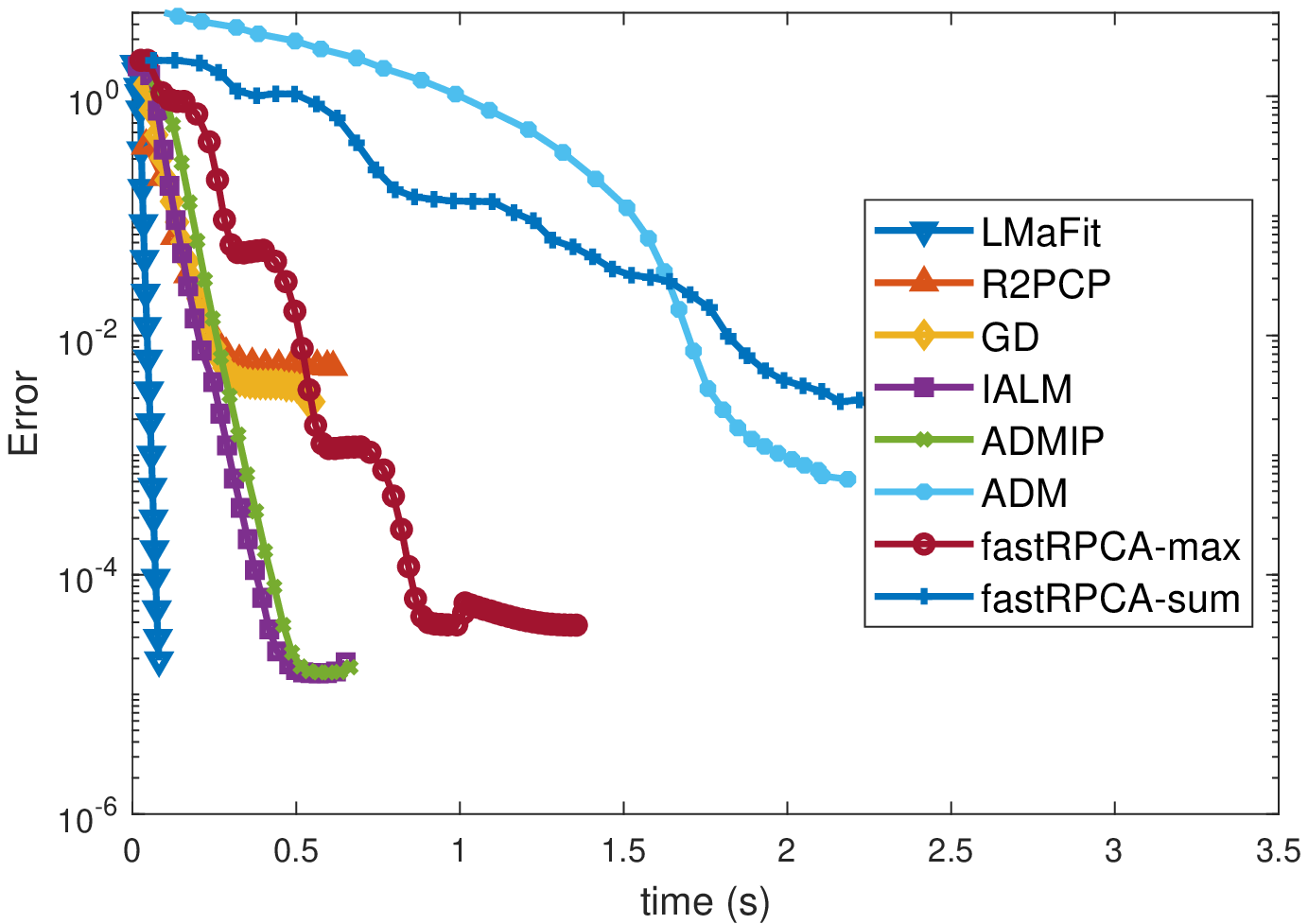}
  \includegraphics[width=0.24\textwidth]{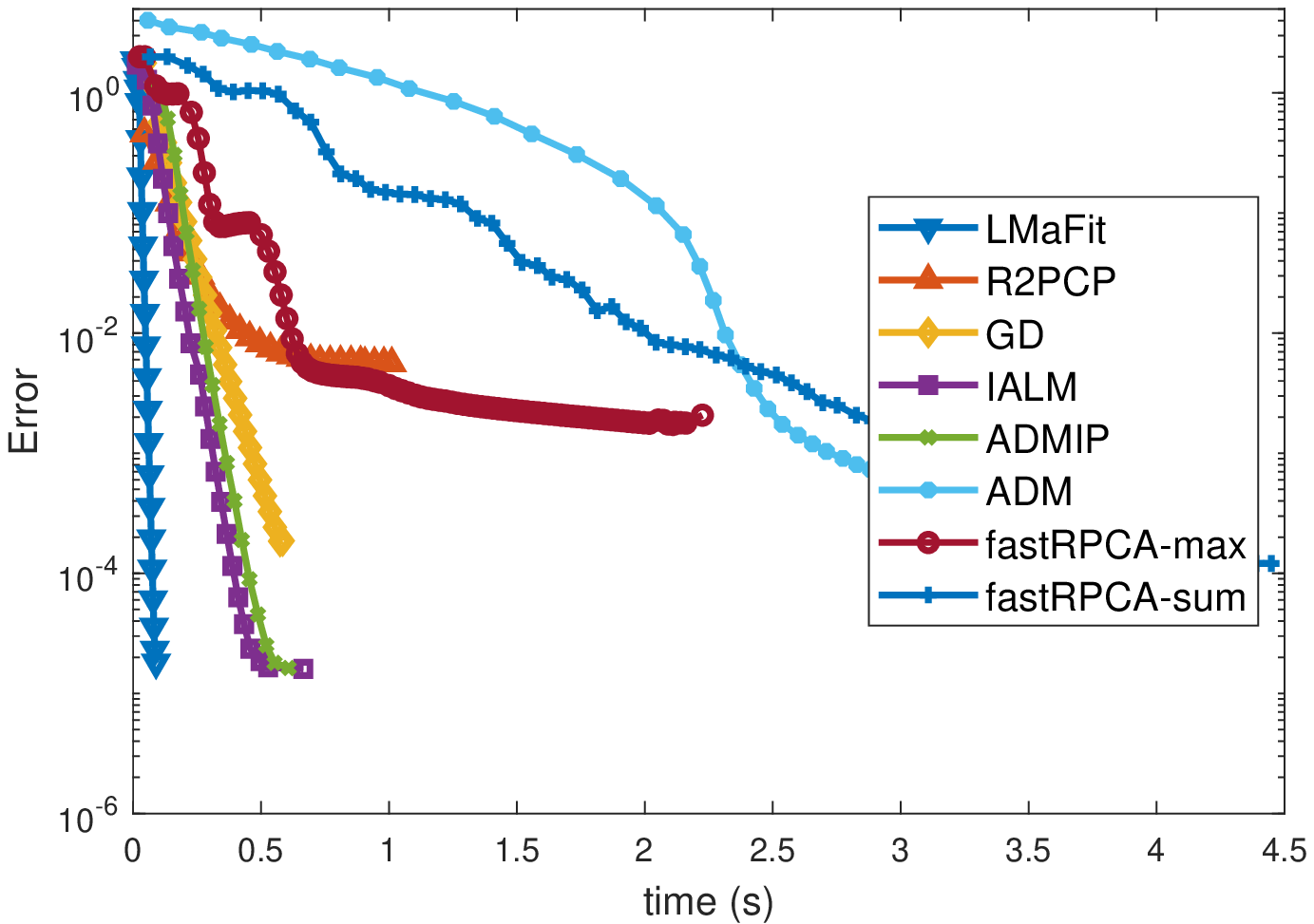}
  \includegraphics[width=0.24\textwidth]{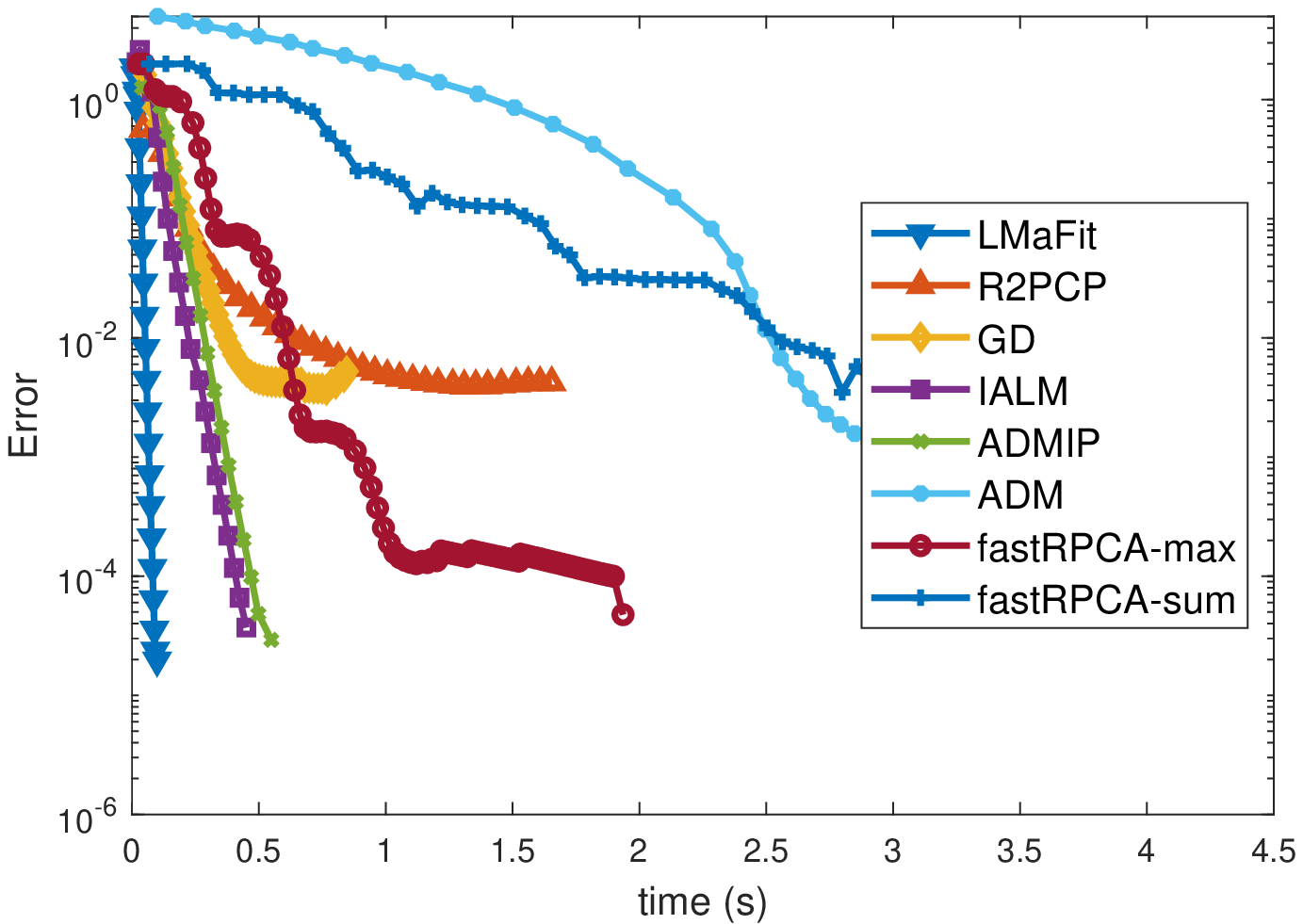}
  \includegraphics[width=0.24\textwidth]{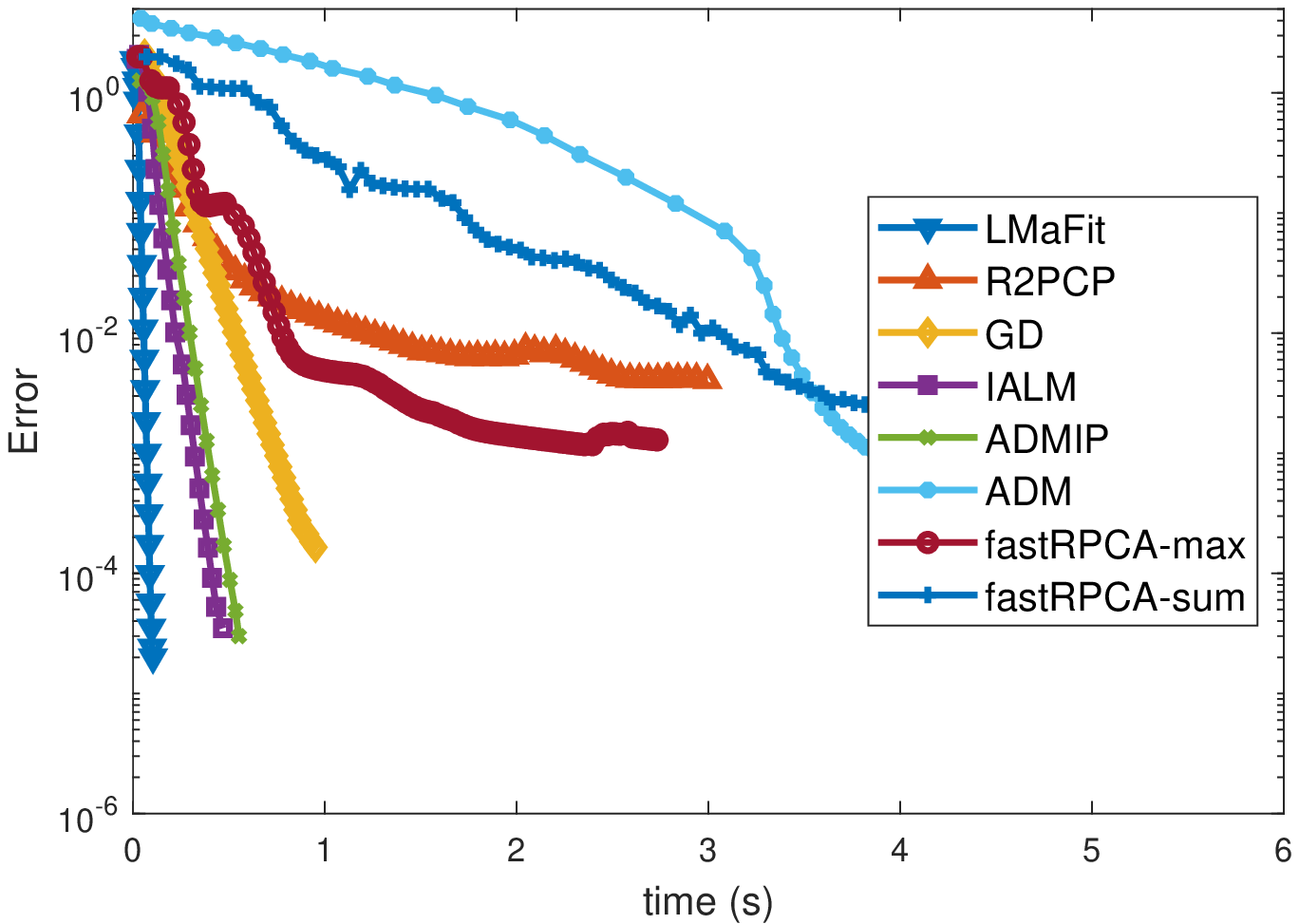}
  \includegraphics[width=0.24\textwidth]{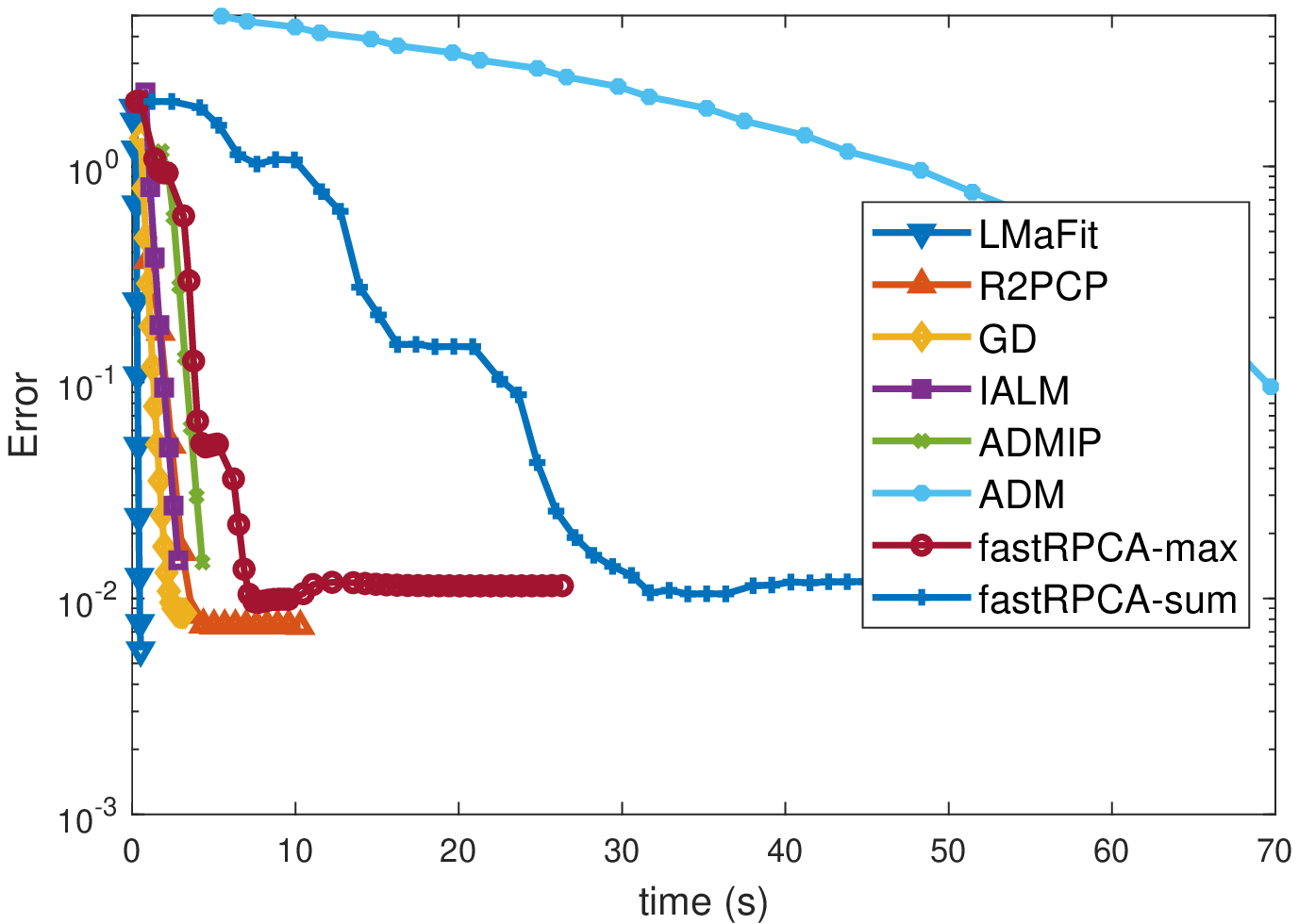}
  \includegraphics[width=0.24\textwidth]{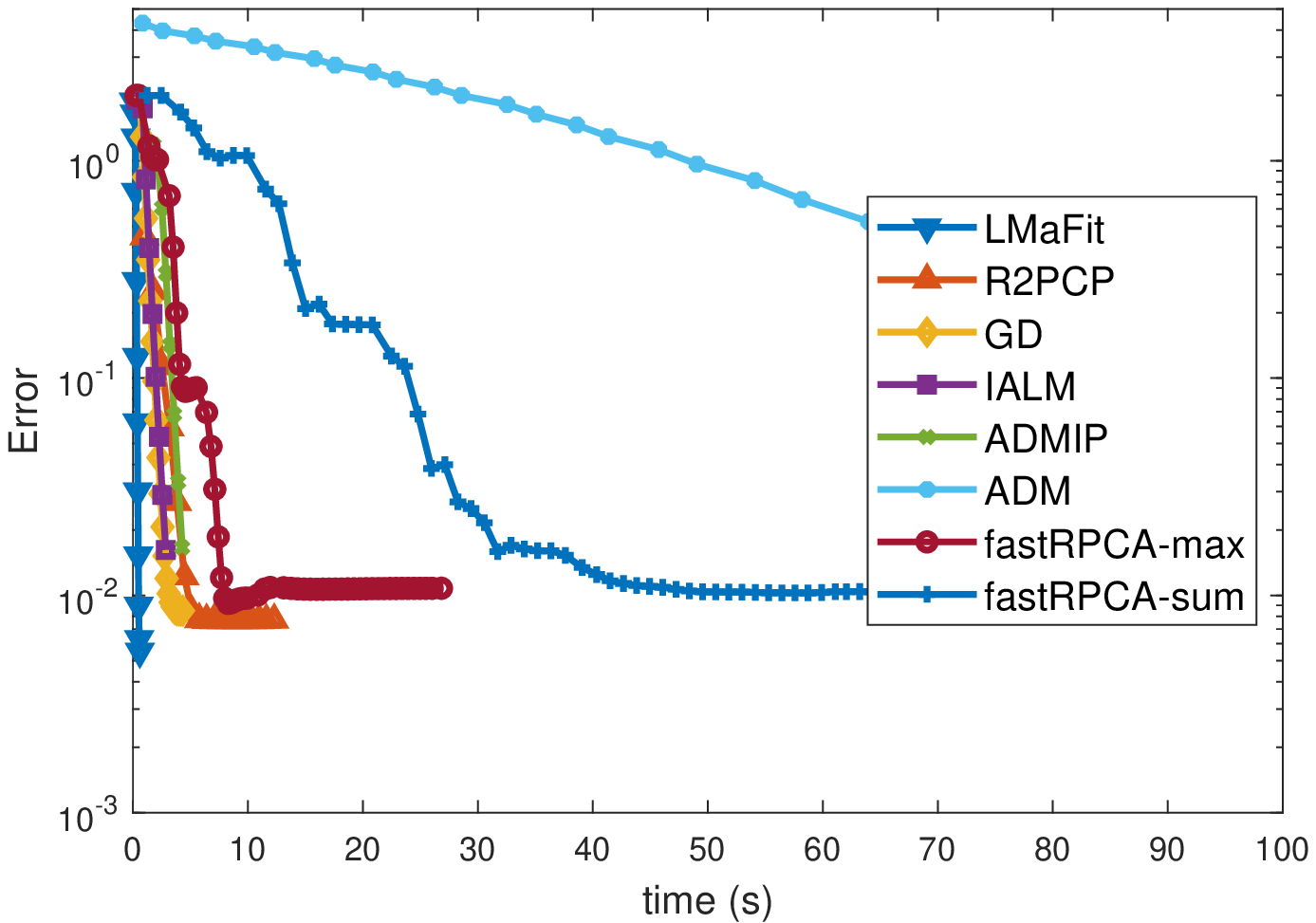}
  \includegraphics[width=0.24\textwidth]{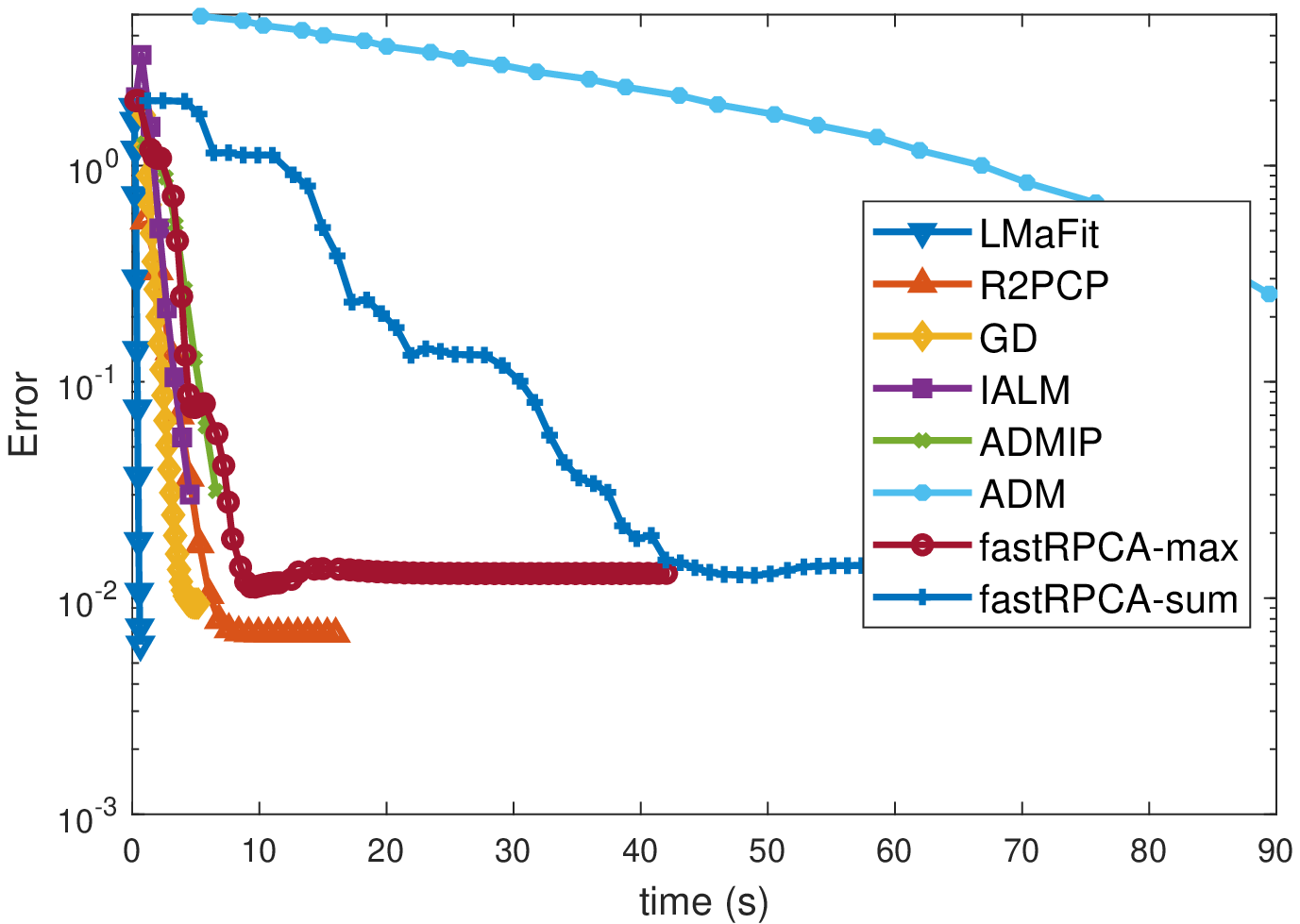}
  \includegraphics[width=0.24\textwidth]{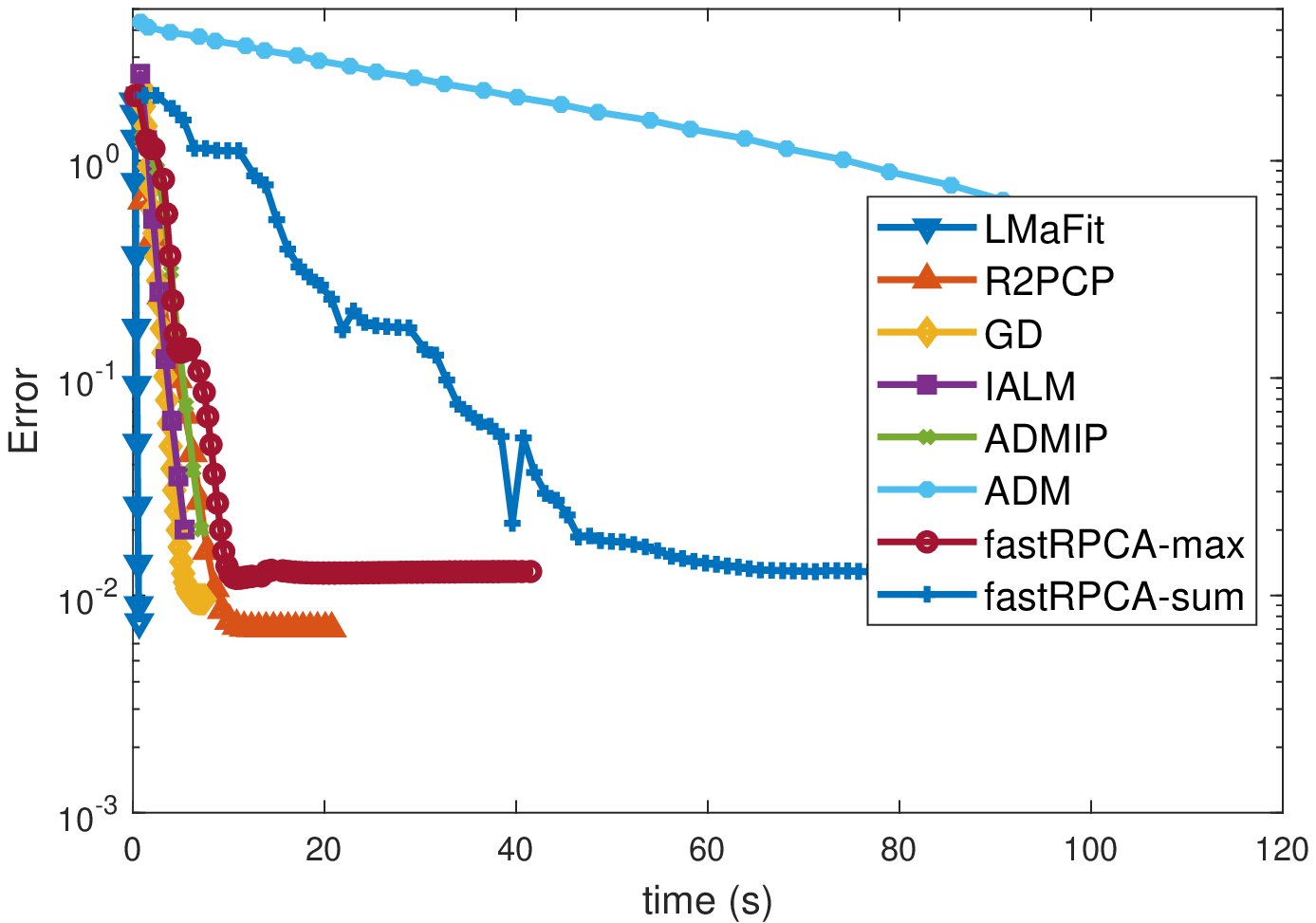}
  \includegraphics[width=0.24\textwidth]{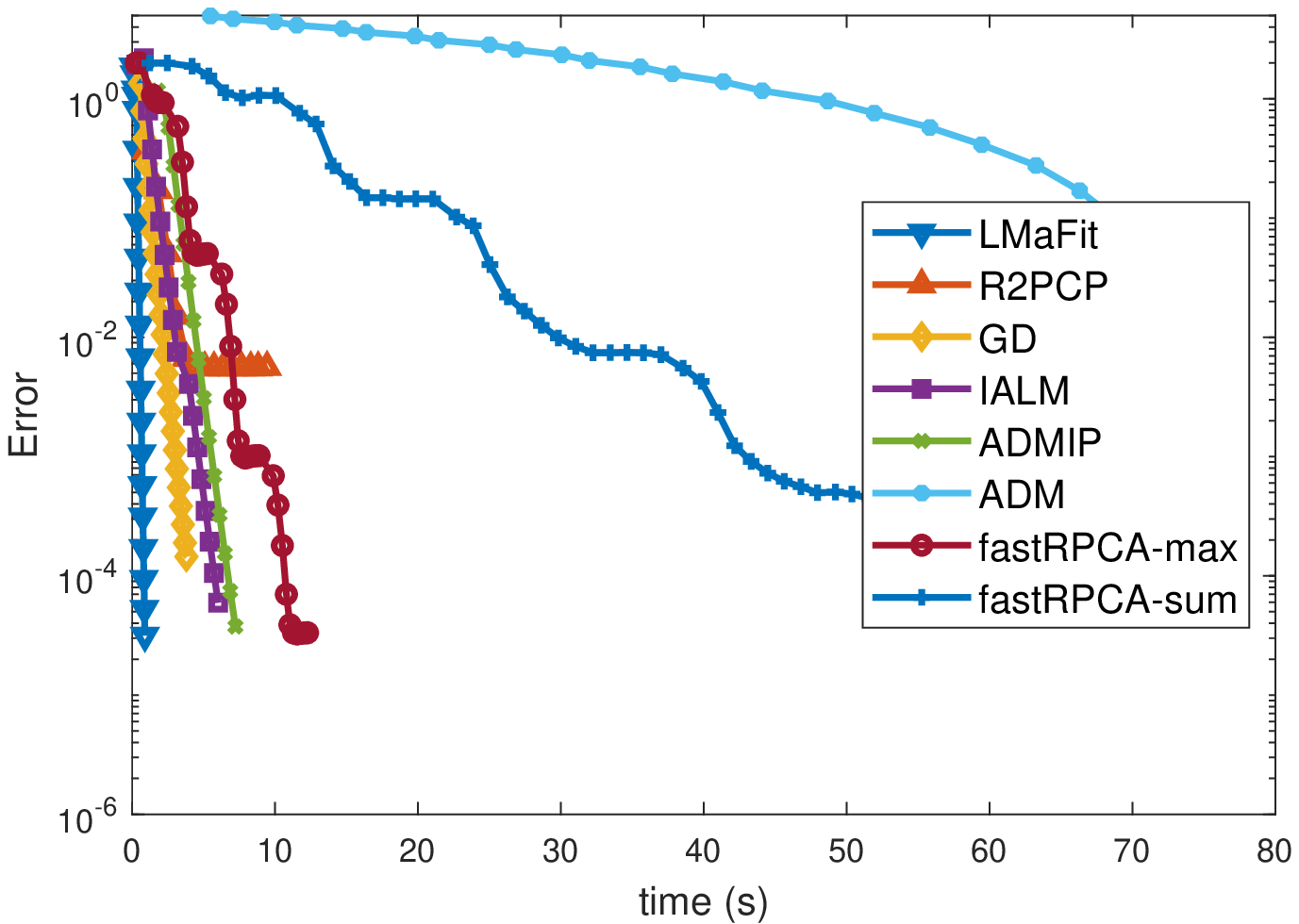}
  \includegraphics[width=0.24\textwidth]{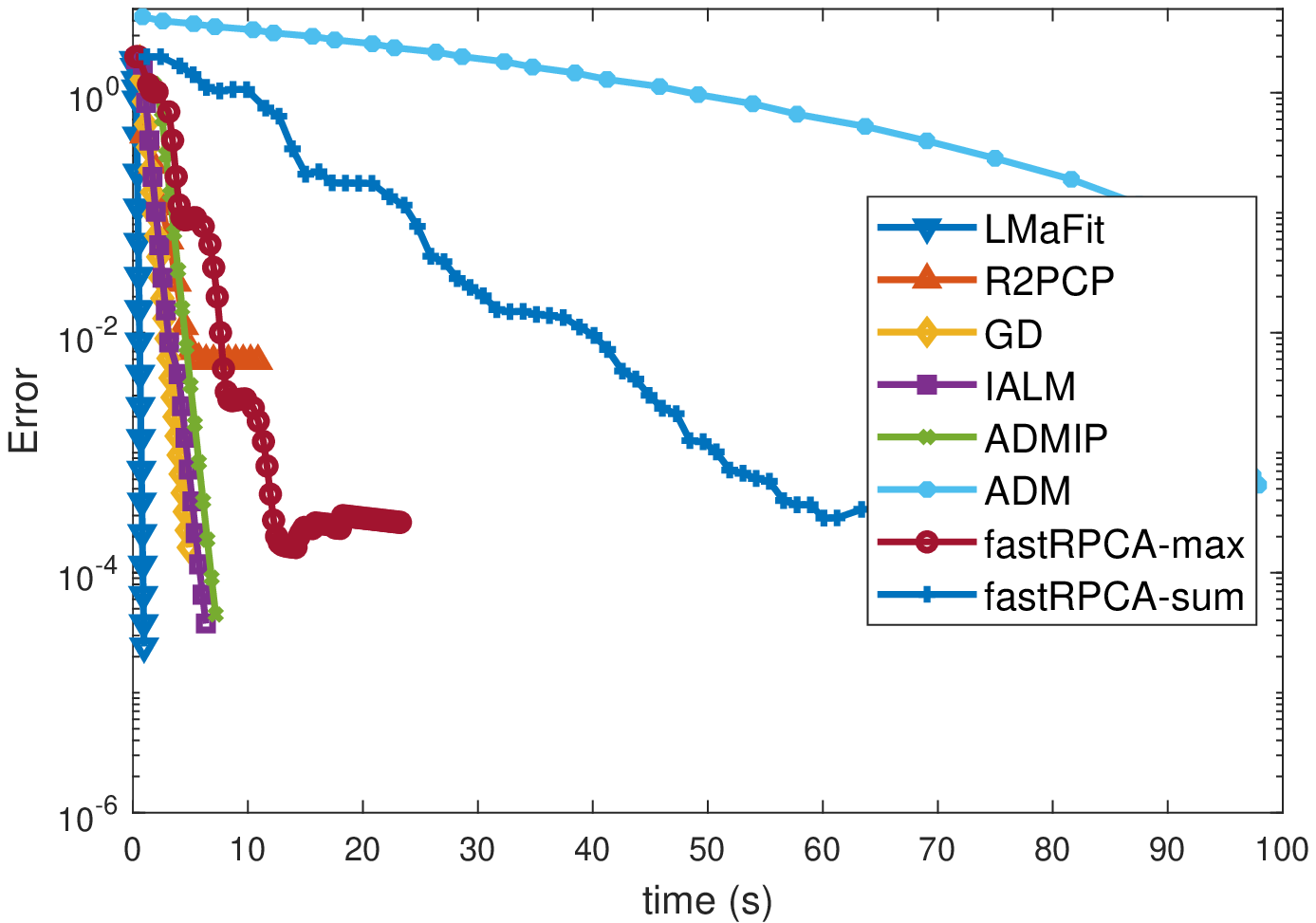}
  \includegraphics[width=0.24\textwidth]{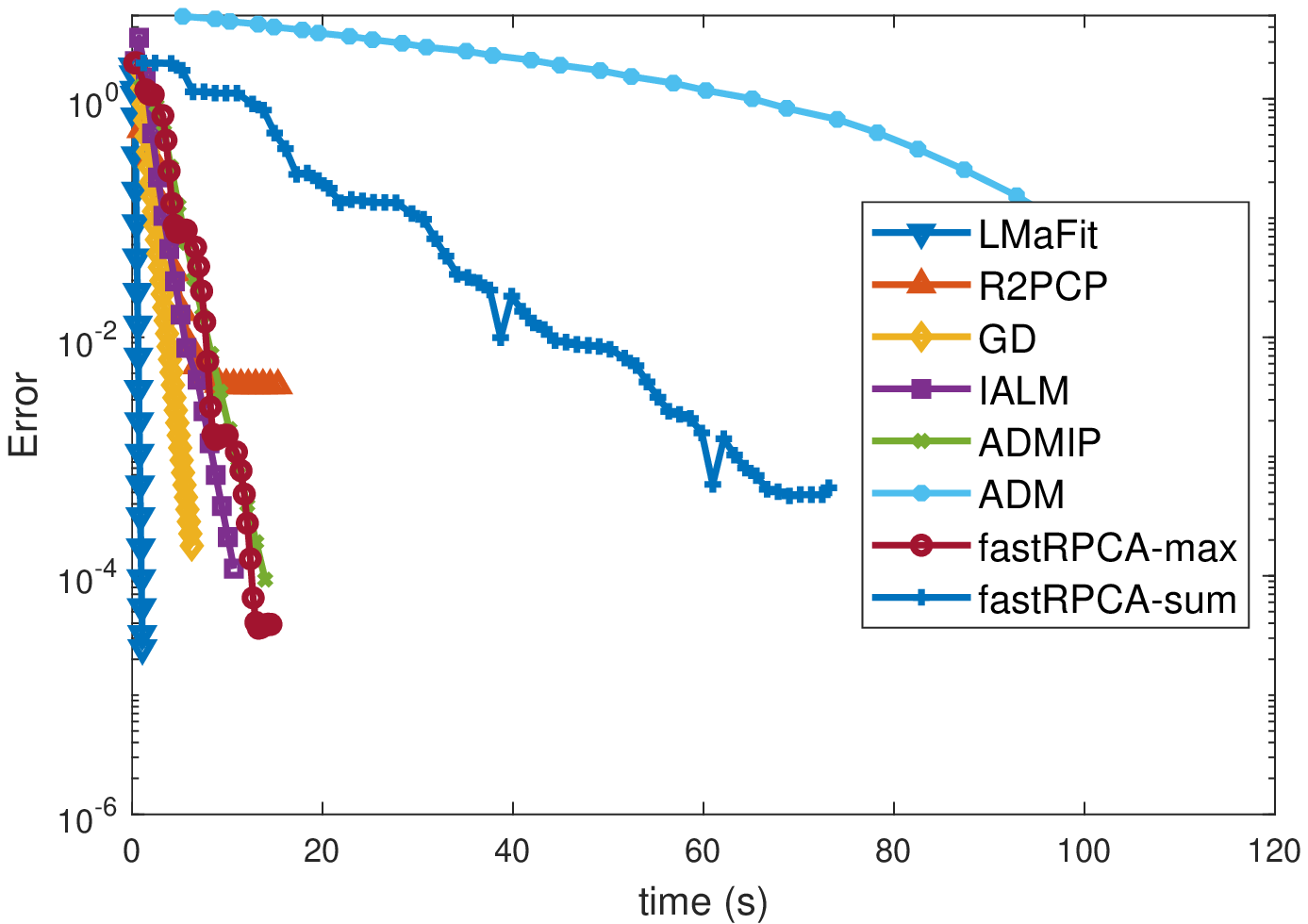}
  \includegraphics[width=0.24\textwidth]{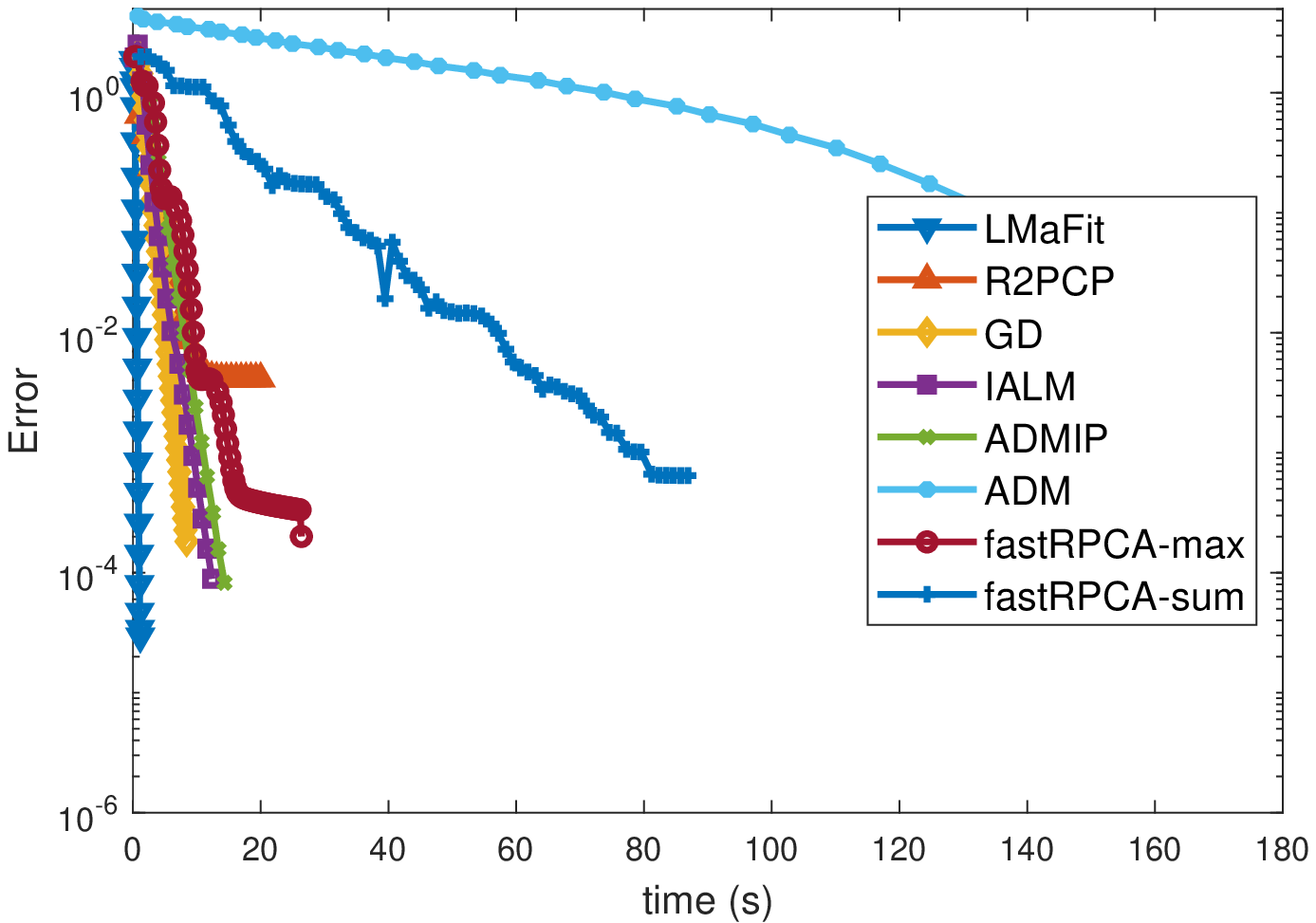}
\end{figure*}

{
In this section we provide some 
{elementary} numerical results of different algorithms for solving RPCA. We selected eight different solvers, 
five for solving convex problems: IALM \cite{Lin-Chen-Wu-Ma-iadm-2009}, ADM \cite{Yuan-Yang-pjo-2013}, ADMIP \cite{Aybat-Iyengar-2015}, fastRPCA-max \cite{P1C2-Aravkin14_1P} and fastRPCA-sum \cite{P1C2-Aravkin14_1P}, and three for solving nonconvex problems: LMafit \cite{Shen-Wen-Zhang-Lmafit}, R2PCP \cite{hintermuller2015robust}, GD \cite{Yi-Park-Chen-Caramanis-gradient}. We tested their performance on some standard synthetic data used in many RPCA papers. The synthetic data were generated in the following manner:
\begin{enumerate}
\item $L^o=UV^\top$ such that $U\in\mR^{n\times r}$, $V\in\mR^{n\times r}$ for $r=c_rn$ and $c_r\in\{0.05,0.1\}$. Moreover, $U_{ij}\sim\cN(0,1)$, $V_{ij}\sim\cN(0,1)$ for all $i,j$ are independent standard Gaussian variables,
\item $\Omega\subset\{(i,j):1\leq i,j\leq n\}$ was chosen uniformly at random such that its cardinality $|\Omega|=c_pn^2$ and $c_p\in\{0.05,0.1\}$,
\item ${S^o_{ij}}\sim\cU[-\sqrt{8r/\pi},\sqrt{8r/\pi}]$ for all $(i,j)\in\Omega$ are independent uniform random variables,
\item $N_{ij}^o\sim\varrho\cN(0,1)$ for all $i,j$ are independent Gaussian variables, where for given signal-to-noise ratio (SNR) of $M$, $\varrho$ is computed from
    \begin{align*}SNR(M) & = 10\log_{10}\left(\frac{\mE[\|L^o+S^o\|_F^2]}{\mE[\|N^o\|_F^2]}\right) \\ & = 10\log_{10}\left(\frac{c_rn+c_s\frac{8r}{3\pi}}{\varrho^2}\right), \end{align*}
    and $SNR(M)\in\{50dB,100dB\}$,
\item The data matrix $M=L^o+S^o+N^o$.
\end{enumerate}
{Note that the non-zero entries of the sparse component and the entries of the low-rank component have approximately the same magnitude in expectation. Indeed, for $n\gg 1$, $L^o_{ij}\approx\sqrt{r}\cN(0,1)$; hence, $\mathbb{E}[L^o_{ij}]=\sqrt{\frac{2r}{\pi}}$ for all $i,j$ and $\mathbb{E}[S^o_{ij}]=\sqrt{\frac{2r}{\pi}}$ for $(i,j)\in\Omega$.}

We created 10 random problems of size $n\in\{500, 1500\}$, i.e., $M\in\mR^{n\times n}$, for each of the two choices of $SNR(M)$, $c_r$ and $c_p$ using the procedure described above. We plot the figures showing the averaged relative errors of the {iterates} over 10 runs versus cpu times in Figure \ref{fig:synthetic}, where the relative error of {$(L,S)$} is defined as
{\[{\rm Error}(L,S) := \frac{\|L-L^o\|_F}{\|L^o\|_F} + \frac{\|S-S^o\|_F}{\|S^o\|_F}.\]}
For all the eight algorithms, we used their default stopping criteria and default parameters setting if the output is of good quality; {otherwise, we tuned some parameters so that the algorithm becomes competitive for our experimental setting. It is worth emphasizing that fastRPCA-sum and fastRPCA-max solve \eqref{eq:aravkin_value_function} with $\psi=\psi_{\rm sum}$ and $\psi=\psi_{\max}$, respectively; and $\rho_{\max}$, the trade-off parameter for $\psi_{\max}$, is set to $\norm{L^o}_*/\norm{S^o}_1$, i.e., this model needs an oracle that provides an ideal separation to tune the parameter.}
\begin{remark}\label{rem:comparison}
We remark that comparing different algorithms for solving RPCA is not an easy task for the following reasons. (i) The algorithms are designed for solving {related but different formulations.} 
For example, IALM solves \eqref{rpca-nuclear-L1} and LMafit solves \eqref{Lmafit-problem-rewrite}, so it is difficult to compare which algorithm is better. (ii) The performance of all the algorithms depends on the parameter settings, initial points, and data structures. For example, from Figure \ref{fig:synthetic} we see that LMafit consistently outperforms fastRPCA-sum. However, this is only based on the data and parameters we tested. For other data sets, initial points, and parameter settings, fastRPCA-sum may be better than LMafit.
\end{remark}
}

\section{Future Directions}
Although \eqref{rpca-nuclear-L1} is a convex problem with nice statistical properties, $\norm{\cdot}_*$ causes problem for large-scale problems as standard algorithms for \eqref{rpca-nuclear-L1} cannot easily exploit multi-core or multi-server computing environments as \eqref{rpca-nuclear-L1} is not amenable to distributed optimization. Some future work in this direction are related to the non-convex formulation which uses an equivalent representation of $\norm{\cdot}_*$ by Recht et al.~\cite{recht2010guaranteed} as given in~\eqref{equiv-nuc-nonconvex}. For instance, given some $\nu>0$, consider $g_\nu:\reals^{m\times n}\rightarrow\reals$ defined in~\eqref{eq:smooth-g}.
Given a data matrix $M\in\reals^{m\times n}$, using $g_\nu$, we can formulate a smooth non-convex optimization problem:
\BE
\label{eq:non-convex-smooth-formulation}
\BA{ll}
\displaystyle\min_{U\in\reals^{m\times r}, V\in\reals^{n\times r}} \ \Psi_\nu(U,V)\triangleq & \tfrac{1}{2}\norm{U}_F^2+\tfrac{1}{2}\norm{V}_F^2 \\ & +\rho g_{\nu}(UV^\top-M),
\EA
\EE
where $\rho>0$ and $\mathbb{Z}_+\ni r\geq \Rank(L^\circ)$ are given parameters.
Here, one can use PALM algorithm~\cite{bolte2014proximal} to generate a sequence that converges to a critical point of $\Psi_\nu$ which is a KL function \sa{-- also see \cite{Markopoulos-TSP-2014,Jiang-Lin-Ma-Zhang-nonconvex-2016} for some other related work on nonconvex optimization.}

\sa{Note that instead of solving the smooth approximation given in \eqref{eq:non-convex-smooth-formulation}, it is preferable to solve the following nonconvex formulation in~\eqref{eq:non-convex-formulation}, which is equivalent to \eqref{rpca-nuclear-L1}.}
\begin{align}\label{eq:non-convex-formulation}
\min_{U\in\reals^{m\times r}, V\in\reals^{n\times r}} &\ \tfrac{1}{2}\norm{U}_F^2+\tfrac{1}{2}\norm{V}_F^2+\rho\norm{UV^\top-M}_1.
\end{align}
\sa{To the best of authors' knowledge, there does not exist efficient methods with strong convergence guarantees to solve \eqref{eq:non-convex-formulation}.}

Note that the third term in \eqref{eq:non-convex-formulation} is a composite function of the form $g(h(\cdot))$ where $g$ is a nonsmooth convex function such that $g$ is Lipschitz continuous, and $h$ is a differentiable function such that its Jacobian $h'$ is Lipschitz. Thus, one possible direction is to design trust region algorithm for \eqref{eq:non-convex-formulation} -- See Section 7.7 in~\cite{ruszczynski2006nonlinear}.
As an alternative to trust-region algorithm, one might also consider the augmented Lagrangian~(AL) method. \red{It is known that for constrained non-convex problems, provided that the second-order KKT conditions hold, the AL will have a saddle point for penalty parameter chosen sufficiently large; therefore, the duality gap encountered in Lagrangian formulations do not pose a problem for augmented Lagrangian based methods -- thus, AL methods might prove useful to establish convergence to local minima when initialized sufficiently close to the local minimum~\cite{ruszczynski2006nonlinear,cohen1984decomposition}.}


{Due to the close relationship between low-rank matrix completion and RPCA, some algorithms for solving low-rank matrix completion problems can possibly be extended to solve variants of RPCA. New methods based on manifold optimization are recently studied for solving low-rank matrix completion problems, see e.g., \cite{OptSpace-2009,Balzano-2010,RTRMC-2011,Vandereycken_2013,Wei-Riemannian-MC-2016}.
It is noted that all these works consider a matrix completion variant/reformulation which is a manifold optimization problem with a smooth objective function. For example, assume that the matrix $M$ is observed partially, i.e., only entries that in a subset $\Omega$ are observed, the low-rank matrix completion model considered in \cite{RTRMC-2011} is:
\BE\label{RTRMC-prob}
\min_{U\in\cG(m,r),V\in\mR^{n\times r}}\half\sum_{(i,j)\in\Omega}C_{ij}^2((UV^\top)_{ij}-M_{ij})^2+\frac{\mu^2}{2}(UV^\top)_{ij}^2,
\EE
where $\cG(m,r)$ denotes Grassmann manifold, $C_{ij}$ denotes some weighting parameter and $\mu>0$ is a penalty parameter. It is noted that \eqref{RTRMC-prob} is a manifold optimization problem with a smooth objective function. Many existing algorithms can be used to solve a manifold optimization problem with smooth objective, for example, Riemannian gradient method \cite{Absil-book}, Riemannian trust-region method \cite{RTRMC-2011} and Riemannian conjugate gradient method \cite{Wei-Riemannian-MC-2016}, and so on.

However, it is more challenging to design algorithms for manifold optimization reformulations of RPCA variants. The reason is that RPCA variants always involve nonsmooth terms in the objective. In fact, all RPCA variants we discussed so far involve $\|S\|_1$ in the objective. As a result, any manifold optimization reformulation of RPCA variants will involve the nonsmooth $\ell_1$ term $\|S\|_1$ as well, unless one can bear with smoothing it, which can potentially degrade the sparsity of $S$. Algorithms for solving manifold optimization problem with nonsmooth objective function have been very limited, and most of them lack convergence guarantees. Nonetheless, some of these algorithms can still be adopted to solve manifold optimization reformulations of RPCA variants, although their efficiency in practice needs further investigations. \sa{For instance, references \cite{hintermuller2015robust,Podosinnikova} propose optimization methods on matrix manifolds. In particular, in \cite{Podosinnikova}, Podosinnikova, Setzer and Hein proposed a new RPCA model by minimizing the trimmed reconstruction error, which reduces to minimizing a nonsmooth function over the Stiefel manifold. The method lacks theoretical convergence guarantees such as convergence to a critical point. That said, the authors of \cite{Podosinnikova} numerically demonstrate that their method exhibits good empirical recovery and it is competitive against other nonconvex formulations and convex optimization based methods. }

One simple manifold optimization reformulation of RPCA is given as follows.
\BE\label{manifold-nonsmooth}
\min_{U\in\St(m,r),V\in\mR^{n\times r}} \ \|UV^\top-M\|_1,
\EE
where $\St(m,r)$ denotes Stiefel manifold. The advantages of \eqref{manifold-nonsmooth} are as follows: (i) it does not involve nuclear norm and thus avoids SVD; (ii) the sizes of $U$ and $V$ are $m\times r$ and $n\times r$, respectively, which are much smaller than the size of $M$ when $r\ll \min(m,n)$. One may also note that \eqref{manifold-nonsmooth} differs \eqref{Lmafit-problem} only with the Stiefel manifold constraint. The drawback of \eqref{Lmafit-problem} is that its optimal solution $(U^*,V^*)$ is not unique, because $(U^*W, V^*W^{-\top})$ is also optimal for any invertible matrix $W\in\mR^{r\times r}$. This drawback is fixed nicely in \eqref{manifold-nonsmooth}. There are several ways to solve \eqref{manifold-nonsmooth}. For example, one can reformulate \eqref{manifold-nonsmooth} as the following one and then apply ADMM to solve it.
\BE\label{manifold-nonsmooth-reformulate}
\BA{ll}
\min & \|S\|_1 \\
\st, & S + UV^\top = M, U\in\St(m,r).
\EA
\EE
The ADMM {iterates the updates} as follows.
\BE\label{manifold-nonsmooth-reformulate-admm}
\BA{ll}
U^{k+1} & := \argmin_U \ \cL_\beta(U,V^k,S^k;\Lambda^k), \st, U\in\St(m,r) \\
V^{k+1} & := \argmin_V \ \cL_\beta(U^{k+1},V,S^k;\Lambda^k) \\
S^{k+1} & := \argmin_S \ \cL_\beta(U^{k+1},V^{k+1},S;\Lambda^k) \\
\Lambda^{k+1} & := \Lambda^k - \beta(U^{k+1}{V^{k+1}}^\top+S^{k+1}-M),
\EA
\EE
where the augmented Lagrangian function is defined as
\[\cL_\beta(U,V,S;\Lambda):=\|S\|_1-\langle UV^\top+S-M\rangle+\frac{\beta}{2}\|UV^\top+S-M\|_F^2.\]
The $U$-subproblem in \eqref{manifold-nonsmooth-reformulate-admm} is a smooth manifold optimization problem and can be solved by existing methods \cite{Absil-book}. This method should be very efficient, but the main issue is that under what kind of conditions it is guaranteed to converge.

Zhang, Ma and Zhang studied some ADMM variants for Riemannian manifold optimization in \cite{Zhang-Ma-Zhang-manifold-2017}, which can be used to solve manifold optimization reformulations of some RPCA variants. We here briefly discuss this work. We consider the following RPCA variant, which minimizes a nonsmooth function over Stiefel manifold.
\BE\label{manifold-problem-admm-ZMZ}
\BA{ll}
\min & \half\|L-UV^\top\|_F^2 + \rho\|S\|_1 + \frac{\mu}{2}\|N\|_F^2 \\
\st  & L + S + N = M \\
     & U \in \St(m,r),
\EA
\EE
where $\rho>0$, $\mu>0$ are tradeoff parameters. One of the ADMM variants for solving \eqref{manifold-problem-admm-ZMZ} proposed in \cite{Zhang-Ma-Zhang-manifold-2017} {iterates the updates} as follows.
\BE\label{manifold-problem-admm-ZMZ-alg}
\BA{ll}
L^{k+1} := & \argmin_L \ \tilde{\cL}_{L^k}(L,U^{k},V^{k},S^{k},N^{k};\Lambda^k), \\ & \st, U \in \St(m,r) \\
U^{k+1} := & \argmin_U \ \tilde{\cL}_{U^k}(L^{k+1},U,V^{k},S^{k},N^{k};\Lambda^k) \\
V^{k+1} := & \argmin_V \ \tilde{\cL}_{V^k}(L^{k+1},U^{k+1},V,S^{k},N^{k};\Lambda^k) \\
S^{k+1} := & \argmin_S \ \tilde{\cL}_{S^k}(L^{k+1},U^{k+1},V^{k+1},S,N^{k};\Lambda^k) \\
N^{k+1} := & N^k - {\eta}\nabla_N\cL(L^{k+1},U^{k+1},V^{k+1},S^{k+1},N;\Lambda^k) \\
\Lambda^{k+1} &:= \Lambda^k - \beta(L^{k+1}+S^{k+1}+N^{k+1}-M),
\EA
\EE
where ${\eta}>0$ is a step size, the augmented Lagrangian function $\cL$ is defined as
\begin{eqnarray*}
\lefteqn{\cL(L,U,V,S,N;\Lambda)  := }\\
& & \half\|L-UV^\top\|_F^2 + \rho_1\|S\|_1 \\
& & + \frac{\mu}{2}\|N\|_F^2  -\langle\Lambda,L + S + N - M \rangle \\
& & + \frac{\beta}{2}\|L + S + N - M\|_F^2,
\end{eqnarray*}
and $\tilde{\cL}$ denotes $\cL$ plus a proximal term. For example, $\tilde{\cL}_{L^k}$ is defined as
\begin{eqnarray*}
\lefteqn{\tilde{\cL}_{L^k}(L,U^{k},V^{k},S^{k},N^{k};\Lambda^k) :=}\\
& & \cL(L,U^{k},V^{k},S^{k},N^{k};\Lambda^k) + \half\|L-L^k\|_H^2,
\end{eqnarray*}
where $H$ denotes a pre-specified positive definite matrix which needs to satisfy certain conditions to guarantee the convergence of \eqref{manifold-problem-admm-ZMZ-alg}. Zhang, Ma and Zhang showed in \cite{Zhang-Ma-Zhang-manifold-2017} that the algorithm described in \eqref{manifold-problem-admm-ZMZ-alg} finds an $\epsilon$-stationary solution to \eqref{manifold-problem-admm-ZMZ} in no more than $O(1/\epsilon^2)$ iterations under certain conditions on $\beta$, $\eta$ and $H$.

One thing that we need to note is that the term $\frac{\mu}{2}\|N\|_F^2$ in \eqref{manifold-problem-admm-ZMZ} is very crucial here. Without this squared term, the convergence results in \cite{Zhang-Ma-Zhang-manifold-2017} do not apply. For example, if one considers the following RPCA variant without the noisy term $N$,
\BE\label{manifold-problem-admm-ZMZ-no-N}
\BA{ll}
\min & \half\|L-UV^\top\|_F^2 + \rho\|S\|_1 \\
\st  & L + S = M \\
     & U \in \St(m,r),
\EA
\EE
then the ADMM variants proposed in \cite{Zhang-Ma-Zhang-manifold-2017} are not guaranteed to converge if they are applied to solve \eqref{manifold-problem-admm-ZMZ-no-N}.
How to extend and generalize the results in \cite{Zhang-Ma-Zhang-manifold-2017} so that they can be used to solve other manifold optimization reformulations of RPCA variants definitely deserves more investigations.
}

{
It is known that the nuclear norm minimization problem can be equivalently written as {an SDP -- 
see \cite{recht2010guaranteed}.} Though the problem size of the SDP is larger than the original nuclear norm minimization problem, it is recently found that the facial reduction technique \cite{Drusvyatskiy2016} can reduce the size of the SDP significantly. This idea has been explored in low-rank matrix completion \cite{Huang-2017} and RPCA \cite{Ma-FR-2018}. In particular, in \cite{Ma-FR-2018} the authors showed that RPCA with partial observation
\BE\label{RPCA-FR}
\min \ \rank(L) + \rho\|S\|_0, \ \st, \ P_{\Omega}(L+S-M) = 0,
\EE
is equivalent to another optimization problem with semidefinite constraint. By applying the facial reduction technique to the semidefinite cone, the size of this reformulation can be significantly reduced, and then it can be solved very efficiently to high accuracy. Extending the facial reduction technique to other RPCA variants is an interesting topic for future research.
}

\section{Conclusions}

In this {paper}, we gave a comprehensive review on algorithms for solving relaxations and variants of robust PCA. {Algorithms for solving convex and nonconvex models were discussed.} We elaborated in details on the applicability of the algorithms and their convergence behaviors. We also proposed several new directions in the hope that they may shed some light for future research in this area.

\section*{Acknowledgments}

The authors are grateful to three anonymous referees for their insightful and constructive comments that led to an improved version of the paper.
The work of S. Ma was supported in part by a startup package from Department of Mathematics at UC Davis. The work of N. S. Aybat was supported in part by NSF under Grant CMMI-1400217 and Grant CMMI-1635106, and ARO grant W911NF-17-1-0298.

\bibliographystyle{IEEEtran}
\bibliography{paper}

\begin{thebibliography}{10}
\providecommand{\url}[1]{#1}
\csname url@samestyle\endcsname
\providecommand{\newblock}{\relax}
\providecommand{\bibinfo}[2]{#2}
\providecommand{\BIBentrySTDinterwordspacing}{\spaceskip=0pt\relax}
\providecommand{\BIBentryALTinterwordstretchfactor}{4}
\providecommand{\BIBentryALTinterwordspacing}{\spaceskip=\fontdimen2\font plus
\BIBentryALTinterwordstretchfactor\fontdimen3\font minus
  \fontdimen4\font\relax}
\providecommand{\BIBforeignlanguage}[2]{{%
\expandafter\ifx\csname l@#1\endcsname\relax
\typeout{** WARNING: IEEEtran.bst: No hyphenation pattern has been}%
\typeout{** loaded for the language `#1'. Using the pattern for}%
\typeout{** the default language instead.}%
\else
\language=\csname l@#1\endcsname
\fi
#2}}
\providecommand{\BIBdecl}{\relax}
\BIBdecl

\bibitem{Candes-Li-Ma-Wright-RPCA}
E.~J. Cand\`{e}s, X.~Li, Y.~Ma, and J.~Wright, ``Robust principal component
  analysis?'' \emph{Journal of ACM}, vol.~58, no.~1, pp. 1--37, 2009.

\bibitem{Peng-rasl-2012}
Y.~Peng, A.~Ganesh, J.~Wright, W.~Xu, and Y.~Ma, ``Rasl: Robust alignment by
  sparse and low-rank decomposition for linearly correlated images,''
  \emph{IEEE Transactions on Pattern Analysis and Machine Intelligence},
  vol.~34, no.~11, pp. 2233--2246, 2012.

\bibitem{Liu-robust-subspace-recovery-2013}
G.~Liu, Z.~Lin, S.~Yan, J.~Sun, Y.~Yu, and Y.~Ma, ``Robust recovery of subspace
  structures by low-rank representation,'' \emph{IEEE Transactions on Pattern
  Analysis and Machine Intelligence}, vol.~35, no.~1, pp. 171--184, 2013.

\bibitem{Shahid-iccv}
N.~Shahid, V.~Kalofolias, X.~Bresson, M.~Bronstein, and P.~Vandergheynst,
  ``Robust principal component analysis on graphs,'' in \emph{ICCV}, 2015.

\bibitem{Chandrasekaran-RPCA-2011}
V.~Chandrasekaran, S.~Sanghavi, P.~A. Parrilo, and A.~S. Willsky,
  ``Rank-sparsity incoherence for matrix decomposition,'' \emph{SIAM Journal on
  Optimization}, vol.~21, no.~2, pp. 572--596, 2011.

\bibitem{NIPS2009_3704}
\BIBentryALTinterwordspacing
J.~Wright, A.~Ganesh, S.~Rao, Y.~Peng, and Y.~Ma, ``Robust principal component
  analysis: Exact recovery of corrupted low-rank matrices via convex
  optimization,'' in \emph{Advances in Neural Information Processing Systems
  22}, Y.~Bengio, D.~Schuurmans, J.~D. Lafferty, C.~K.~I. Williams, and
  A.~Culotta, Eds.\hskip 1em plus 0.5em minus 0.4em\relax Curran Associates,
  Inc., 2009, pp. 2080--2088. [Online]. Available:
  \url{http://papers.nips.cc/paper/3704-robust-principal-component-analysis-exact-recovery-of-corrupted-low-rank-matrices-via-convex-optimization.pdf}
\BIBentrySTDinterwordspacing

\bibitem{NIPS2010_4005}
\BIBentryALTinterwordspacing
H.~Xu, C.~Caramanis, and S.~Sanghavi, ``Robust {PCA} via outlier pursuit,'' in
  \emph{Advances in Neural Information Processing Systems 23}, J.~D. Lafferty,
  C.~K.~I. Williams, J.~Shawe-Taylor, R.~S. Zemel, and A.~Culotta, Eds.\hskip
  1em plus 0.5em minus 0.4em\relax Curran Associates, Inc., 2010, pp.
  2496--2504. [Online]. Available:
  \url{http://papers.nips.cc/paper/4005-robust-pca-via-outlier-pursuit.pdf}
\BIBentrySTDinterwordspacing

\bibitem{recht2010guaranteed}
B.~Recht, M.~Fazel, and P.~A. Parrilo, ``Guaranteed minimum-rank solutions of
  linear matrix equations via nuclear norm minimization,'' \emph{SIAM review},
  vol.~52, no.~3, pp. 471--501, 2010.

\bibitem{Zhou-stable-pca-2010}
Z.~Zhou, J.~Wright, X.~Li, E.~J. Candès, and Y.~Ma, ``Stable principal
  component pursuit,'' in \emph{Proceedings of International Symposium on
  Information Theory,}, 2010.

\bibitem{Zhang-zhou-Liang}
H.~Zhang, Y.~Zhou, and Y.~Liang, ``Analysis of robust {PCA} via local
  incoherence,'' 2015.

\bibitem{Yi-Park-Chen-Caramanis-gradient}
X.~Yi, D.~Park, Y.~Chen, and C.~Caramanis, ``Fast algorithms for robust {PCA}
  via gradient descent,'' 2016.

\bibitem{Zhang-Yang-manifold-2017}
T.~Zhang and Y.~Yang, ``Robust {PCA} by manifold optimization,''
  \emph{https://arxiv.org/abs/1708.00257}, 2017.

\bibitem{Netrapalli-nonconvex-2014}
P.~Netrapalli, U.~N. Niranjan, S.~Sanghavi, A.~Anandkumar, and P.~Jain,
  ``Non-convex robust {PCA},'' 2014.

\bibitem{Bouwmans-survey-cviu-2014}
T.~Bouwmans and E.~H. Zahzah, ``Robust {PCA} via principal component pursuit: A
  review for a comparative evaluation in video surveillance,'' \emph{Computer
  Vision and Image Understanding}, vol. 122, pp. 22--34, 2014.

\bibitem{Aybat-book-2016}
N.~S. Aybat, \emph{Handbook of Robust Low Rank and Sparse Matrix Decomposition:
  Applications in Image and Video Processing}.\hskip 1em plus 0.5em minus
  0.4em\relax CRC Press, Taylor and Francis Group, 2016, ch. Algorithms for
  Stable PCA.

\bibitem{Lin-WCAMSAP-2009}
Z.~Lin, A.~Ganesh, J.~Wright, L.~Wu, M.~Chen, and Y.~Ma, ``Fast convex
  optimization algorithms for exact recovery of a corrupted low-rank matrix,''
  in \emph{International Workshop on Computational Advances in Multi-Sensor
  Adaptive Processing}, 2009.

\bibitem{Lin-Chen-Wu-Ma-iadm-2009}
Z.~Lin, M.~Chen, L.~Wu, and Y.~Ma, ``The augmented lagrange multiplier method
  for exact recovery of corrupted low-rank matrices,'' UIUC Technical Report
  UILU-ENG-09-2215, Tech. Rep., 2009.

\bibitem{Yuan-Yang-pjo-2013}
X.~Yuan and J.~Yang, ``Sparse and low-rank matrix decomposition via alternating
  direction methods,'' \emph{Pacific Journal of Optimization}, vol.~9, no.~1,
  pp. 167--180, 2013.

\bibitem{Goldfarb-Ma-Scheinberg-2013}
D.~Goldfarb, S.~Ma, and K.~Scheinberg, ``Fast alternating linearization methods
  for minimizing the sum of two convex functions,'' \emph{Mathematical
  Programming Series A}, vol. 141, no. 1-2, pp. 349--382, 2013.

\bibitem{Tao-Yuan-siopt-2011}
M.~Tao and X.~Yuan, ``Recovering low-rank and sparse components of matrices
  from incomplete and noisy observations,'' \emph{SIAM Journal on
  Optimization}, vol.~21, no.~1, pp. 57--81, 2011.

\bibitem{Aybat-Ma-Goldfarb-2014}
N.~S. Aybat, D.~Goldfarb, and S.~Ma, ``Efficient algorithms for robust and
  stable principal component pursuit problems,'' \emph{Computational
  Optimization and Applications}, vol.~58, pp. 1--29, 2014.

\bibitem{Aybat-Iyengar-2015}
N.~S. Aybat and G.~Iyengar, ``An alternating direction method with increasing
  penalty for stable principal component pursuit,'' \emph{Computational
  Optimization and Applications}, vol.~61, pp. 635--668, 2015.

\bibitem{P1C2-Aravkin14_1P}
A.~Aravkin, S.~Becker, V.~Cevher, and P.~Olsen, ``A variational approach to
  stable principal component pursuit,'' in \emph{30th Conference on Uncertainty
  in Artificial Intelligence (UAI) 2014}, no. EPFL-CONF-199542, 2014, preprint
  made available at arXiv:1406.1089 [math.OC] in June, 2014.

\bibitem{Lin-Ma-Zhang-2015-free-gamma}
T.~Lin, S.~Ma, and S.~Zhang, ``Global convergence of unmodified 3-block {ADMM}
  for a class of convex minimization problems,'' \emph{Journal of Scientific
  Computing}, 2017.

\bibitem{Mu-FW-2016}
C.~Mu, Y.~Zhang, J.~Wright, and D.~Goldfarb, ``Scalable robust matrix recovery:
  Frank-wolfe meets proximal methods,'' \emph{SIAM Journal on Scientific
  Computing}, vol.~38, no.~5, pp. A3291--A3317, 2016.

\bibitem{GoDec-icml-2011}
T.~Zhou and D.~Tao, ``Godec: Randomized low-rank and sparse matrix
  decomposition in noisy case,'' in \emph{ICML}, 2011.

\bibitem{pmlr-v31-zhou13b}
------, ``Greedy bilateral sketch, completion \& smoothing,'' in
  \emph{Proceedings of the Sixteenth International Conference on Artificial
  Intelligence and Statistics}, ser. Proceedings of Machine Learning Research,
  C.~M. Carvalho and P.~Ravikumar, Eds., vol.~31.\hskip 1em plus 0.5em minus
  0.4em\relax Scottsdale, Arizona, USA: PMLR, 29 Apr--01 May 2013, pp.
  650--658.

\bibitem{hintermuller2015robust}
M.~Hinterm{\"u}ller and T.~Wu, ``Robust principal component pursuit via inexact
  alternating minimization on matrix manifolds,'' \emph{Journal of Mathematical
  Imaging and Vision}, vol.~51, no.~3, pp. 361--377, 2015.

\bibitem{Gu-Wang-Liu-AISTATS-2016}
Q.~Gu, Z.~Wang, and H.~Liu, ``Low-rank and sparse structure pursuit via
  alternating minimization,'' in \emph{AISTATS}, 2016.

\bibitem{Xu-OR-PCA-2013}
J.~Feng, H.~Xu, and S.~Yan, ``Online robust {PCA} via stochastic
  optimization,'' in \emph{NIPS}, 2013.

\bibitem{Shen-Wen-Zhang-Lmafit}
Y.~Shen, Z.~Wen, and Y.~Zhang, ``Augmented lagrangian alternating direction
  method for matrix separation based on low-rank factorization,''
  \emph{Optimization Methods and Software}, vol.~29, no.~2, pp. 239--263, 2014.

\bibitem{Jiang-Lin-Ma-Zhang-nonconvex-2016}
B.~Jiang, T.~Lin, S.~Ma, and S.~Zhang, ``Structured nonconvex and nonsmooth
  optimization: Algorithms and iteration complexity analysis.''
  \emph{https://arxiv.org/abs/1605.02408}, 2016.

\bibitem{Nesterov-07}
Y.~E. Nesterov, ``Gradient methods for minimizing composite functions,''
  \emph{Mathematical Programming}, vol. 140, no.~1, pp. 125--161, 2013.

\bibitem{Beck-Teboulle-2009}
A.~Beck and M.~Teboulle, ``A fast iterative shrinkage-thresholding algorithm
  for linear inverse problems,'' \emph{SIAM J. Imaging Sciences}, vol.~2,
  no.~1, pp. 183--202, 2009.

\bibitem{Tseng-2008}
P.~Tseng, ``On accelerated proximal gradient methods for convex-concave
  optimization,'' \emph{Manuscript}, 2008.

\bibitem{Combettes-2007}
P.~L. Combettes and J.-C. Pesquet, ``A douglas-rachford splitting approach to
  nonsmooth convex variational signal recovery,'' \emph{IEEE Journal of
  Selected Topics in Signal Processing}, vol.~1, no.~4, pp. 564--574, 2007.

\bibitem{Goldstein-Osher-split-Bregman-2009}
T.~Goldstein and S.~Osher, ``The split {B}regman method for l1-regularized
  problems,'' \emph{SIAM Journal on Imaging Sciences}, vol.~2, no.~2, pp.
  323--343, 2009.

\bibitem{Yang-Zhang-Yin-ADMM-2010}
J.~Yang, Y.~Zhang, and W.~Yin, ``A fast alternating direction method for
  tvl1-l2 signal reconstruction from partial fourier data,'' \emph{IEEE Journal
  of Selected Topics in Signal Processing Special Issue on Compressed Sensing},
  vol.~4, no.~2, pp. 288--297, 2010.

\bibitem{Yang-Zhang-admm-2011}
J.~Yang and Y.~Zhang, ``Alternating direction algorithms for l1-problems in
  compressive sensing,'' \emph{SIAM Journal on Scientific Computing}, vol.~33,
  no.~1, pp. 250--278, 2011.

\bibitem{Gabay-Mercier-1976}
D.~Gabay and B.~Mercier, ``A dual algorithm for the solution of nonlinear
  variational problems via finite-element approximations,'' \emph{Comp. Math.
  Appl.}, vol.~2, pp. 17--40, 1976.

\bibitem{Gabay-83}
D.~Gabay, ``Applications of the method of multipliers to variational
  inequalities,'' in \emph{Augmented Lagrangian Methods: Applications to the
  Solution of Boundary Value Problems}, M.~Fortin and R.~Glowinski, Eds.\hskip
  1em plus 0.5em minus 0.4em\relax Amsterdam: North-Hollan, 1983.

\bibitem{Fortin-Glowinski-1983}
M.~Fortin and R.~Glowinski, \emph{Augmented Lagrangian methods: applications to
  the numerical solution of boundary-value problems}.\hskip 1em plus 0.5em
  minus 0.4em\relax North-Holland Pub. Co., 1983.

\bibitem{Lions-Mercier-79}
P.~L. Lions and B.~Mercier, ``Splitting algorithms for the sum of two nonlinear
  operators,'' \emph{SIAM Journal on Numerical Analysis}, vol.~16, pp.
  964--979, 1979.

\bibitem{Eckstein-Bertsekas-1992}
J.~Eckstein and D.~P. Bertsekas, ``On the {D}ouglas-{R}achford splitting method
  and the proximal point algorithm for maximal monotone operators,''
  \emph{Mathematical Programming}, vol.~55, pp. 293--318, 1992.

\bibitem{He-Yuan-rate-ADM-2012}
B.~He and X.~Yuan, ``On the $\mathcal{O}(1/n)$ convergence rate of
  douglas-rachford alternating direction method,'' \emph{SIAM Journal on
  Numerical Analysis}, vol.~50, pp. 700--709, 2012.

\bibitem{Monteiro-Svaiter-2010a}
R.~D.~C. Monteiro and B.~F. Svaiter, ``Iteration-complexity of
  block-decomposition algorithms and the alternating direction method of
  multipliers,'' \emph{SIAM Journal on Optimization}, vol.~23, pp. 475--507,
  2013.

\bibitem{Lin-Ma-Zhang-convergence-2014}
T.~Lin, S.~Ma, and S.~Zhang, ``On the sublinear convergence rate of multi-block
  {ADMM},'' \emph{Journal of the Operations Research Society of China}, vol.~3,
  no.~3, pp. 251--274, 2015.

\bibitem{Chen-admm-failure-2013}
C.~Chen, B.~He, Y.~Ye, and X.~Yuan, ``The direct extension of {ADMM} for
  multi-block convex minimization problems is not necessarily convergent,''
  \emph{Mathematical Programming}, vol. 155, pp. 57--79, 2016.

\bibitem{aybat2011fast}
N.~S. Aybat, D.~Goldfarb, and G.~Iyengar, ``Fast first-order methods for stable
  principal component pursuit,'' \emph{arXiv preprint arXiv:1105.2126}, 2011.

\bibitem{larsen1998lanczos}
R.~M. Larsen, ``Lanczos bidiagonalization with partial reorthogonalization,''
  \emph{DAIMI Report Series}, vol.~27, no. 537, 1998.

\bibitem{P1C2-Aravkin13_1J}
A.~Y. Aravkin, J.~V. Burke, and M.~P. Friedlander, ``Variational properties of
  value functions,'' \emph{SIAM Journal on optimization}, vol.~23, no.~3, pp.
  1689--1717, 2013.

\bibitem{FW-1956}
M.~Frank and P.~Wolfe, ``An algorithm for quadratic programming,'' \emph{Naval
  Res. Logis. Quart.}, vol.~3, pp. 95--110, 1956.

\bibitem{Jaggi-icml-2013}
M.~Jaggi, ``Revisiting frank-wolfe: Projection-free sparse convex
  optimization,'' in \emph{ICML}, 2013.

\bibitem{thesisMJ}
------, ``Sparse convex optimization methods for machine learning,'' Ph.D.
  dissertation, ETH Zurich, Oct. 2011.

\bibitem{Li-Ng-Yuan-NLAA-2015}
X.~Li, M.~K. Ng, and X.~Yuan, ``Median filtering-based methods for static
  background extraction from surveillance video,'' \emph{Numerical Linear
  Algebra with Applications}, vol.~22, pp. 845--865, 2015.

\bibitem{Li-Wang-Hu-Cai-2014}
L.~Li, P.~Wang, Q.~Hu, and S.~Cai, ``Efficient background modeling based on
  sparse representation and outlier iterative removal,'' \emph{IEEE
  Transactions on Circuits and Systems for Video Technology}, vol.~26, no.~2,
  pp. 278--289, 2014.

\bibitem{Ma-Sparcoc-2015}
S.~Ma, D.~Johnson, C.~Ashby, D.~Xiong, C.~L. Cramer, J.~H. Moore, S.~Zhang, and
  X.~Huang, ``{SPARCoC}: a new framework for molecular pattern discovery and
  cancer gene identification,'' \emph{PLoS ONE}, vol.~10, no.~3, p. e0117135,
  2015.

\bibitem{Yang-Pong-Chen-2017}
L.~Yang, T.~K. Pong, and X.~Chen, ``Alternating direction method of multipliers
  for a class of nonconvex and nonsmooth problems with applications to
  background/foreground extraction,'' \emph{SIAM J. Imaging Sciences}, vol.~10,
  pp. 74--110, 2017.

\bibitem{Kurdyka-1998}
K.~Kurdyka, ``On gradients of functions definable in o-minimal structures,''
  \emph{Annales de l'institut Fourier}, vol. 146, pp. 769--783, 1998.

\bibitem{Lojasiewicz-1963}
S.~{\L}ojasiewicz, \emph{Une propri$\acute{e}$t$\acute{e}$ topologique des
  sous-ensembles analytiques r$\acute{e}$els, Les $\acute{E}$quations aux
  D$\acute{e}$riv$\acute{e}$es Partielles}.\hskip 1em plus 0.5em minus
  0.4em\relax Paris: $\acute{E}$ditions du centre National de la Recherche
  Scientifique, 1963.

\bibitem{Lewis-Malick-alt-proj-manifold}
A.~S. Lewis and J.~Malick, ``Alternating projections on manifolds,''
  \emph{Mathematics of Operations Research}, vol.~33, no.~1, pp. 216--234,
  2008.

\bibitem{sparcs-2011}
A.~E. Waters, A.~C. Sankaranarayanan, and R.~Baraniuk, ``Sparcs: Recovering
  lowrank and sparse matrices from compressive measurements,'' in \emph{NIPS},
  2011.

\bibitem{glz13}
S.~Ghadimi, G.~Lan, and H.~Zhang, ``Mini-batch stochastic approximation methods
  for nonconvex stochastic composite optimization,'' \emph{Math. Program.},
  vol. 155, no.~1, pp. 267--305, 2016.

\bibitem{bolte2014proximal}
J.~Bolte, S.~Sabach, and M.~Teboulle, ``Proximal alternating linearized
  minimization for nonconvex and nonsmooth problems,'' \emph{Mathematical
  Programming}, vol. 146, no. 1-2, pp. 459--494, 2014.

\bibitem{Markopoulos-TSP-2014}
P.~P. Markopoulos, G.~N. Karystinos, and D.~A. Pados, ``Optimal algorithms for
  {L}1-subspace signal processing,'' \emph{IEEE Transactions on Signal
  Processing}, vol.~62, pp. 5046--5058, 2014.

\bibitem{ruszczynski2006nonlinear}
A.~P. Ruszczy{\'n}ski, \emph{Nonlinear optimization}.\hskip 1em plus 0.5em
  minus 0.4em\relax Princeton university press, 2006, vol.~13.

\bibitem{cohen1984decomposition}
G.~Cohen and D.~L. Zhu, ``Decomposition coordination methods in large scale
  optimization problems. the nondifferentiable case and the use of augmented
  lagrangians,'' \emph{Advances in large scale systems}, vol.~1, pp. 203--266,
  1984.

\bibitem{OptSpace-2009}
R.~H. Keshavan and S.~Oh, ``Optspace: A gradient descent algorithm on the
  grassman manifold for matrix completion,'' \emph{Arxiv preprint
  arXiv:0910.5260v2}, 2009.

\bibitem{Balzano-2010}
L.~Balzano, R.~Nowak, and B.~Recht, ``Online identification and tracking of
  subspaces from highly incomplete information,'' in \emph{Allerton
  Conference}, 2010.

\bibitem{RTRMC-2011}
N.~Boumal and P.-A. Absil, ``{RTRMC}: A {R}iemannian trust-region method for
  low-rank matrix completion,'' in \emph{NIPS}, 2011.

\bibitem{Vandereycken_2013}
B.~Vandereycken, ``Low-rank matrix completion by {R}iemannian optimization,''
  \emph{SIAM J. Optim.}, vol.~23, no.~2, pp. 1214--1236, 2013.

\bibitem{Wei-Riemannian-MC-2016}
K.~Wei, J.-F. Cai, T.~F. Chan, and S.~Leung, ``Guarantees of riemannian
  optimization for low rank matrix completion,'' \emph{preprint
  https://arxiv.org/abs/1603.06610}, 2016.

\bibitem{Absil-book}
P.-A. Absil, R.~Mahony, and R.~Sepulchre, \emph{Optimization algorithms on
  matrix manifolds}.\hskip 1em plus 0.5em minus 0.4em\relax Princeton
  University Press, 2008.

\bibitem{Podosinnikova}
A.~Podosinnikova, S.~Setzer, and M.~Hein, ``Robust {PCA}: Optimization of the
  robust reconstruction error over the stiefel manifold,'' \emph{GCPR}, 2014.

\bibitem{Zhang-Ma-Zhang-manifold-2017}
J.~Zhang, S.~Ma, and S.~Zhang, ``Multi-block optimization over {R}iemannian
  manifolds: an iteration complexity analysis,''
  \emph{https://arxiv.org/abs/1710.02236}, 2017.

\bibitem{Drusvyatskiy2016}
D.~Drusvyatskiy and H.~Wolkowicz, ``The many faces of degeneracy in conic
  optimization,'' \emph{Foundations and Trends in Optimization}, vol.~3, no.~2,
  pp. 77--170, 2016.

\bibitem{Huang-2017}
S.~Huang and H.~Wolkowicz, ``Low-rank matrix completion using nuclear norm with
  facial reduction,'' \emph{Journal of Global Optimization}, 2017.

\bibitem{Ma-FR-2018}
S.~Ma, F.~Wang, L.~Wei, and H.~Wolkowicz, ``Robust principal component analysis
  using facial reduction,''
  \emph{http://www.optimization-online.org/DB\_HTML/2018/03/6535.html}, 2018.

\end{thebibliography}

\end{document}